\documentclass[a4paper,12pt]{article} 
\newdimen\paperhight
\usepackage{amsmath}
\usepackage{amssymb}
\usepackage{amsfonts}
\usepackage[usenames]{color}

\newcommand{\ch}{{\rm ch}} 

\newcommand{\pr}{\par \vspace{3mm}\noindent [{\bf Proof}] \qquad}
\newcommand{\prend}{\hfill \qed  \vspace{3mm}}

\newcommand{\qed}{\quad\hbox{\rule[-2pt]{3pt}{6pt}}\par\vspace{3mm}}

\newcommand{\1}{{\bf 1}} 
 
\newcommand{\C}{\mathbb C} 
\newcommand{\Z}{\mathbb Z} 
\newcommand{\Q}{\mathbb Q} 
 
\newcommand{\N}{\mathbb N}

\newcommand{\CA}{{\cal A}}

\newcommand{\CD}{{\cal D}}

\newcommand{\CI}{{\cal I}}

\newcommand{\CO}{{\cal O}}
\newcommand{\CP}{{\cal P}}
\newcommand{\CQ}{{\cal Q}}
\newcommand{\CR}{{\cal R}}
\newcommand{\CS}{{\cal S}}

\newcommand{\CY}{{\cal Y}}

\newcommand{\CH}{{\cal H}}

\newcommand{\wt}{{\rm wt}}
\newcommand{\Hom}{{\rm Hom}}
\newcommand{\End}{{\rm End}}
\newcommand{\Aut}{{\rm Aut}}

\newcommand{\Tr}{{\rm Tr}}
\newtheorem{thm}{Theorem}
\newtheorem{prn}[thm]{Proposition}
\newtheorem{dfn}[thm]{Definition}

\newtheorem{lmm}[thm]{Lemma}

\begin{document}
\title{A $\Z_3$-orbifold theory of lattice vertex operator 
algebra and $\Z_3$-orbifold constructions}
\author{\begin{tabular}{c}
Masahiko Miyamoto\\
Institute of Mathematics, \\
University of Tsukuba, \\
Tsukuba, 305 Japan \end{tabular}}
\date{}
\maketitle

\begin{abstract}
Let $V$ be a simple VOA of CFT-type satisfying $V'\cong V$ and $G$ 
a finite automorphism group of $V$. We prove that 
if all $V$-modules are completely reducible and 
a fixed point subVOA $V^G$ is $C_2$-cofinite, 
then all $V^G$-modules are completely reducible and 
every simple $V^G$-module appears in a 
$g$-twisted (or ordinary) module of $V$ as a 
$V^G$-submodule for some $g\in G$. 
We also prove that $V_L^{\sigma}$ is $C_2$-cofinite 
for a lattice VOA $V_L$ and 
$\sigma\in \Aut(V_L)$ lifted from a triality automorphism of $L$. 
Using these results, we present two $\Z_3$-orbifold constructions as examples. 
One is the moonshine VOA $V^{\natural}$ and the other is a new 
CFT No.32 in Schellekens' list \cite{S}. 
\end{abstract}  

\section{Introduction}
A concept of vertex operator algebras (shortly VOAs) $V=(V,Y,\1,\omega)$ 
was introduced by Borcherds \cite{B} with a purpose to explain  
the moonshine phenomenon \cite{CN} and then as a stage for studying 
the phenomenon, Frenkel, Lepowsky and Meurman \cite{FLMe} constructed 
the moonshine VOA $V^{\natural}$ 
using a $\Z_2$-orbifold construction 
from the Leech lattice VOA $V_{\Lambda}$. 
A VOA is now understood 
as an algebraic version of a $2$-dim. conformal field theory (shortly CFT). 
Among CFTs, a rational CFT has an  
important meaning because one can determine all 
representations precisely, where 
we interpret "rational" as meaning that 
all $\N$-gradable modules are completely reducible. 
Therefore, it is important to find new VOAs whose modules are all completely 
reducible. One way to construct such candidates 
is an orbifold theory. It is a 
theory of the fixed point subVOA $V^G$ given by a finite 
automorphism group $G$ of $V$. 
For example, if all $\N$-gradable $V$-modules are completely reducible, then 
$V^G$ is expected to have the same property. 
Furthermore, $V$ has a special module called a $g$-twisted 
module for $g\in G$ (see \cite{DLiMa}), which 
is a direct sum of $V^{g}$-modules on 
which $V$ acts as a permutation in a sense. 
So there is a natural question if every simple 
$V^G$-module appears as a submodule of a 
$g$-twisted (or ordinary) $V$-module for some $g\in G$. 
These statement are ambiguously expected to be true \cite{DHVW}.  
One of the main purposes in this paper is to show that this is true under 
some conditions called $C_2$-cofiniteness.  
Throughout this paper, we will only treat a simple VOA 
$V=\oplus_{m=0}^{\infty}V_m$ of CFT-type 
with a nonsingular invariant bilinear form, 
that is, $V_0=\C \1$ and $V$ is isomorphic to its restricted dual 
$V'=\oplus_m \Hom(V_m,\C)$. 

Except for the complete reducibility of modules, 
another important condition for VOAs is $C_2$-cofiniteness. 
This is defined by the condition that a subspace 
$$C_2(V):=<v_{-2}u \mid v,u\in V>_{\C}$$ 
of $V$ has a finite co-dimension in $V$, 
where $v_{-2}$ denotes a coefficient of vertex operator 
$$Y(v,z)=\sum_{m\in \Z} v_mz^{-1-m}\in \End(V)[[z,z^{-1}]]$$
of $v\in V$ at $z$. This assumption was introduced 
by Zhu \cite{Z} as a technical condition to prove 
an ${\rm SL}_2(\Z)$-invariance property of 
the space spanned by trace functions $\Psi_W$ on 
$V$-modules $W$ under the assumption 
that all $V$-modules are completely reducible. 
However, as the author has shown in \cite{Mi1}, this is a natural condition 
from a view point of the representation theory 
(it is equivalent to the nonexistence of $\N$-ungradable modules) and 
is enough to prove some kind of ${\rm SL}_2(\Z)$-invariance 
property, that is, 
the author has introduced pseudo-trace functions and the space spanned by 
trace functions and pseudo-trace functions is invariant under the action of 
$SL_2(\Z)$. Frankly speaking, pseudo-trace functions are not easy to treat. 
For example, the author has recently shown in \cite{Mi2} and \cite{Mi3} that  
if a $C_2$-cofinite VOA $V$ satisfies $V'\cong V$ and an $S$-transformation 
$S(\Psi_V)$ of a trace function $\Psi_V$ on $V$ is a linear combination of 
trace functions, then $V$ satisfies some generalized concept 
of rigidity for nonsemisimple modules. 
He called it semirigid, but this is different from the 
classical means. For example, in the definition of rigidity, for 
a surjection $\phi:W\to V$, we consider an injection $V\to W$. 
However, if $W$ is not completely reducible, then the existence 
of surjection $\phi:W\to V$ does not guarantee the existence of 
injection $V\to W$.

\begin{dfn}
$\mbox{}\quad$ $W\in {\rm mod}(V)$ is called to be {\bf semi-rigid} 
if there are $\widetilde{W}\in {\rm mod}(V)$, $e_W\in \Hom_V(W\boxtimes \widetilde{W},V)$ and   
$e_{\widetilde{W}}\in \Hom_V(\widetilde{W}\boxtimes W,V)$ such that \\
(1) $W, \widetilde{W}$ are $C_1$-cofinite as $V$-modules and 
$W\boxtimes \widetilde{W}$ and $W\boxtimes (\widetilde{W}\boxtimes W)$ 
are $V$-modules, \\ 
(2) there are $Q\in {\rm mod}(V)$ and an embedding $\epsilon:Q\to \widetilde{W}\boxtimes W$ such that 
$e_{\widetilde{W}}\epsilon:Q\to V$ is a covering and \\
(3) 
in the diagram
$$\begin{array}{ccccc}
W\boxtimes Q&\xrightarrow{{\rm id}_W\boxtimes {\epsilon}}
&W\boxtimes (\widetilde{W}\boxtimes W)& \xrightarrow{{\rm id}_W\boxtimes e_{\widetilde{W}} } 
&W\boxtimes V \cr
&&\mbox{}\qquad\downarrow \mu  & & \cr
 &&(W\boxtimes \widetilde{W})\boxtimes W   & \xrightarrow{e_W\boxtimes {\rm id}_W} & V\boxtimes W
\end{array}\eqno{(1.1)}$$ 
$(e_W\boxtimes {\rm id}_W)\mu({\rm id}_W\boxtimes \epsilon)$ is surjective,  
where $\mu:W\boxtimes (\widetilde{W}\boxtimes W) \rightarrow (W\boxtimes \widetilde{W})\boxtimes W$ 
is a natural isomorphism for the associativity of 
products of intertwining operators (see (5.3)). 
\end{dfn}

Under his definition of semirigidity, the fusion products preserves 
the exactness of sequences, (see Proposition 1 (vii)), where 
$S$-transformation is given by $\begin{pmatrix}0&-1\cr 1&0\end{pmatrix}$.  
Good examples without pseudo-trace functions are orbifold models. 
For example, the author has also shown in \cite{Mi3} that 
$S(\Psi_{V^g})$ has no pseudo-trace functions for each 
$g\in \Aut(V)$ of finite order. 
We will also prove this is true for $S(\Psi_{V^G})$ for any 
finite automorphism group $G$ of $V$. \\

\noindent
{\bf Theorem A} \quad {\it 
Let $V$ be a simple VOA of CFT-type. 
Assume $V'\cong V$ and all $V$-modules are completely reducible. 
If $G$ is a finite automorphism group of $V$ and 
a fixed point subVOA $V^G$ is $C_2$-cofinite, then $V^G$ is semirigid 
and every simple $V^G$-module appears 
as a $V^G$-submodule of a 
$g$-twisted (or ordinary) $V$-module for some $g\in G$. 
In particular, if $G$ is solvable, then all $V^G$-modules are completely reducible.}\\

Let's show examples of orbifold models. 
For every positive definite even lattice $L$, there is a VOA $V_L$ associated 
with $L$ called a lattice VOA. As well known, all $V_L$-modules are 
completely reducible \cite{D}.  If $\sigma$ acts on $L$ as $-1$, then all 
$V_L^{\sigma}$-modules are classified in \cite{Ys}. 
This result relies heavily on  
a wonderful result done by Dong and Nagatomo \cite{DNa} about   
the fixed pointed subVOA of free bosonic Fock space. 
Unfortunately, there is no such a result for other automorphisms 
at the present time. 

The main object in this paper is an automorphism $\sigma$ of $L$ 
of order three. 
For a special lattice and an automorphism, there is a classification 
of simple modules \cite{TYa}.  We will treat a general case. \\

\noindent
{\bf Theorem B} \quad {\it 
Let $L$ be a positive definite even lattice and $V_L$ a 
lattice VOA associated with $L$.  
Let $\sigma\in \Aut(L)$ of order three. We use the same notation for 
an automorphism of $V_L$ lifted from $\sigma$. 
Then a fixed point subVOA $V_L^{\sigma}$ is $C_2$-cofinite.}\\

At the end of this paper, as applications of these results,  
we will give two $\Z_3$-orbifold constructions as examples. 
One is a VOA which has the same character as $V^{\natural}$ does and 
the other is a new meromorphic $c=24$ VOA which has the same 
character with No. 32 in Schellekens' list \cite{S}.\\

\noindent
{\bf Theorem C} \quad {\it 
Let $\Lambda$ be a positive definite even unimodular lattice of rank $N$ 
with an automorphism $\sigma$ of order three.  
Let $H=\Lambda^{\sigma}$ be a fixed point sublattice of $\Lambda$ 
and assume that $N-{\rm rank}(H)$ is divisible by $3$. Then 
we are able to construct a VOA $\widetilde{V}$ by 
a $\Z_3$-orbifold construction 
from a lattice VOA $V_{\Lambda}$. }

\section{Preliminary results}
\subsection{$C_2$-cofiniteness}
In this situation, we will assume $C_2$-cofiniteness only. 
All modules in this paper are finitely generated. \vspace{-4mm}\\

\begin{prn}\label{prn:C2}
Let $V$ be a $C_2$-cofinite VOA. Then we have the followings: \\
(i) Every module is $\Z_+$-gradable and weights of modules are in {\rm $\Q$, \cite{Mi1}}. \\
(ii) The number of inequivalent simple modules is finite, {\rm \cite{GN},\cite{DLiMa}}.\\
(iii) Set $V=B+C_2(V)$ for $B$ spanned by homogeneous elements. 
Then for any module $W$ generated from one element $w$ has a spanning set 
$\{v^1_{n_1}....v^k_{n_k}w \mid v^i\in B, \quad  n_1<\cdots <n_k\}$. 
Hence every f.g. $V$-module is $C_n$-cofinite for any $n=1,2,...$, {\rm \cite{Mi1},\cite{Bu},\cite{GN}}.\\
(iv) Every $V$-module has a projective cover. \\
(v)  For any $V$-modules $W$ and $U$, a fusion product $W\boxtimes U$ 
is well-defined as f.g. modules. \\
(vi) If $V\cong V'$, then $V$ is projective if and only 
if all modules are completely reducible. \\
(vii) If $V\cong V'$, then 
for any exact sequence $0\to A\to B\to C\to 0$ and a module $W$, 
$0\to W\boxtimes A\to W\boxtimes B \to W\boxtimes C \to 0$ is still exact.
See {\rm \cite{Mi2}, \cite{Mi3}} for (iv)$\sim$ (vii).
\end{prn}

\subsection{Intertwining operators}
For $V$-modules $A,B,C$, let 
$\CI_{A,B}^C$ be the set of (logarithmic) intertwining operators of 
$A$ from $B$ to $C$. 
Since $V$ is $C_2$-cofinite, $\CY$ satisfies a differential 
equation of regular 
singular points and so there is $K\in \N$ such that $\CY$ has a form
$$\CY(a,z)=\sum_{i=0}^K \sum_{n\in \C} a_{(n,i)}z^{-n-1}\log^iz
\in \Hom(B,C)\{z\}[\log z].$$
We note that the action of $L(0)$ on a module may not be semisimple. 
Let $\wt$ denote the semisimple part of $L(0)$. 
$\CY^{(m)}(a,z)=\sum_{n\in \C} a_{(n,m)}z^{-n-1}$ denotes the 
coefficient at $\log^mz$. 
From the $L(-1)$-derivative property for $\CY$, 
we have two important properties:
$$\begin{array}{rl}
\displaystyle{\CY^{(m)}(a,z)=}&\displaystyle{\frac{1}{m!}(z\frac{d}{dz}-zL(-1))^m\CY^{(0)}(a,z)}, 
\qquad \mbox{ and }\cr
\displaystyle{(i+1)a_{(n,i+1)}b
=}&\displaystyle{-(L(0)\!-\!\wt)a_{(n,i)}b+((L(0)\!-\!\wt)a)_{(n,i)}b+a_{(n,i)}((L(0)\!-\!\wt)b)} 
\end{array} $$
for $b\in B$. In particular, ${\CY}^{(K)}(\ast,z)$ is an ordinary 
intertwining operator (i.e. of formal power series).  
One more important result is that $(L(0)-\wt)W$ is a proper 
submodule for a $V$-module $W\not=0$.  

As Huang has shown in \cite{H1}, for $d'\in D', a\in A, b\in B, c\in C$ and 
intertwining operators 
$\CY_1\in \CI_{A,E}^D$, $\CY_2\in \CI_{B,C}^E$, 
$\CY_3\in \CI_{F,C}^D$ and $\CY_4\in \CI_{A,B}^F$,  
the formal power series (with logarithmic terms) 
$$\langle  d',\CY_1(a,x)\CY_2(b,y)c\rangle \quad\mbox{ and }\quad
\langle  d',\CY_3(\CY_4(a,x-y)b,y)c\rangle $$
are absolutely convergent when $|x|>|y|>0$ and $|y|>|x-y|>0$, respectively, 
and can all be analytically extended to multi-valued analytic functions on 
$$M^2=\{(x,y)\in \C^2 \mid xy(x-y)\not=0 \}.$$
As he did, we are able to lift them to single-valued analytic functions 
$$
E(\langle  d,\CY_1(a,x)\CY_2(b,y)c\rangle ) \quad\mbox{ and }\quad
E(\langle  d,\CY_3(\CY_4(a,x-y)b,y)c\rangle ) $$
on the universal covering $\widetilde{M^2}$ of $M^2$. 
As he remarked, single-valued liftings are not unique, but 
the existence of such functions is enough for our arguments. 
The important fact is that if we fix $A,B,C,D$, then these functions 
are given as solutions of the same differential equations. 
Therefore, for $\CY_1\in \CI_{A,E}^D,\CY_2\in \CI_{B,C}^E$ 
there are $\CY_5\in \CI_{A\boxtimes B,C}^D$ and 
$\CY_6\in \CI_{B,A\boxtimes C}^D$ such that 
$$\begin{array}{cl}
E(\langle d',\CY_1(a,x)\CY_2(b,y)c\rangle)
=E(\langle d',\CY_5(\CY_{A,B}^{\boxtimes}(a,x-y)b,y)c\rangle) &\mbox{ and}\cr
E(\langle d',\CY_1(a,x)\CY_2(b,y)c\rangle)
=E(\langle d',\CY_6(b,y)\CY_{A,C}^{\boxtimes}(a,x)c\rangle), &
\end{array}$$
where $\CY_{A,B}^{\boxtimes}$ denotes an intertwining operator to 
define a fusion product $A\boxtimes B$.

\subsection{VOA whose modules are all completely reducible} 
In this subsection, we will explain very important known properties  
of a simple $C_2$-cofinite VOA $V$ whose modules are all 
completely reducible and $V\cong V'$. Let $N$ be 
the central charge of $V$ and $\{V\cong W^0,\cdots,W^s\}$ the set 
of all simple $V$-modules.\\  

\noindent
{\bf 2.3.1  Zhu's modular invariance property}\\
For $v\in V_m$ with $L(1)v=0$, we consider a trace function 
$$  T_{W^i}(v;\tau)=\Tr_{W^i} o(v)q^{L(0) -N/24}, $$
where $q=e^{2\pi \sqrt{-1}\tau}$ and $o(v)=v_{\wt(v)-1}$ 
is a grade-preserving operator of $v$ on $W^i$. 
In particular, a trace function of the vacuum $\1\in V$ 
$$  T_{W^i}(\1;\tau)=q^{-N/24}\sum_{m} \dim (W^i_m) q^m $$
is called a character of $W^i=\oplus_m W^i_m$ and we denote it by $\ch(W^i)$.
As Zhu has shown in \cite{Z}, these functions are well-defined 
in the upper half plane $\CH=\{\tau\in \C \mid {\rm Im}(\tau)>0\}$.  
A wonderful property of these VOAs is an ${\rm SL}_2(\Z)$-invariance 
property. Namely, there are $s_{ij}\in \C$ which does not depend 
on $v$ such that 
$$(1/\tau)^{\wt(v)}T_{W^i}(v;-1/\tau)=\sum_{j=0}^s s_{ij}T_{W^j}(v;\tau).$$
We will call the transformation in the left side an $S$-transformation 
of $T_{W^i}$ 
and then the matrix $S=(s_{ij})_{i,j=0,\ldots,s}$ an $S$-matrix of $V$. \\

\noindent
{\bf 2.3.2  Dong, Li and Mason's modular invariance property}\\
Dong, Li and Mason \cite{DLiMa} extended the above result to the case 
where they also consider 
an automorphism $\sigma$ of order $n$. For $g,h\in <\sigma>$, 
they introduced a concept of $g$-twisted $h$-stable modules $W$. 
Let's explain a $\sigma$-twisted module briefly. See \cite{DLiMa} 
for the detail. Decompose 
$$V=\oplus_{i=0}^{n-1} V^{(i)}, \quad \mbox{ where } V^{(i)}=
\{v\in V\mid \sigma(v)=e^{2\pi \sqrt{-1} i/n}v\}.$$
Clearly, $V^{(0)}$ is a subVOA. A simple $\sigma$-twisted module $W$ 
has a grading $W=\oplus_{m=0}^{\infty} W_{r+m/n}$ such that 
every $W^{(i)}=\oplus_{k=0}^{\infty} W_{r+i/n+k}$ 
is a simple $V$-module, where $W_t$ denotes the space consisting of 
elements of weight $t$.  
Furthermore, there is $r\in \Z$ such that 
$$v_t(W^{(i)})\subseteq W^{(i+rp)}$$ 
for $v\in V^{(p)}$ and $t\in \Q$, where we consider $W^{(k)}=W^{(h)}$ 
if $k\equiv h \pmod{n}$. 
We will call $W^{(i)}$ a $V^{(0)}$-sector of $W$.  
On the other hand, the definition of $h$-stable modules $W$ is 
$h\circ W\cong W$ as $V$-modules, 
where $v_n(h\circ w)=h\circ (h(v))_n w$ for $w\in W$. 
This implies the existence 
of an endomorphism $\phi(h)$ of $W$ such that 
$\phi(h)^{-1}v_n\phi(h)=(h(v))_n$ for $n\in \Z$ and $v\in V$.  
For such a module, they consider a trace function 
$$  T_W(h,v;\tau)=\Tr_{W} \phi(h)o(v)e^{2\pi \sqrt{-1}(\tau -c/24)} $$
of $\phi(h)o(v)$ on a $g$-twisted $h$-stable module $W$ for $v\in V_m$ 
with $L(1)v=0$.  Then they have shown that for 
$A=\begin{pmatrix}a&b\cr e&f\end{pmatrix}\in {\rm SL}_2(\Z)$, 
$A$-transformation 
$$  (e\tau+f)^{-m} T_{W}(h,v;\frac{a\tau+b}{e\tau+f})$$
of $T_W(h,v;\tau)$ is a linear sum of trace functions of 
$\phi(g^ah^e)o(v)$ on 
$g^ah^e$-stable $g^bh^f$-twisted modules under the assumption that 
all twisted modules are completely reducible.  
In particular, $(\frac{1}{\tau})^mT_V(\sigma,v;-1/\tau)$ is a 
linear sum of trace functions of $o(v)$ 
on $\sigma$-twisted modules. 

An interesting case is that $V$ has exactly one simple module $V$. 
In this case, $V$ has only one simple $\sigma^i$-twisted module 
$V(\sigma^i)$ for each $i$ and 
so they are all $\sigma^j$-stable.  If the weights of elements 
in $V(\sigma^i)$ are in $\Z/n$, then a homomorphism $\phi(\sigma)$ 
for a $\sigma$-stable module is given by 
$\phi(\sigma)=e^{2\pi r\sqrt{-1}L(0)}$ since 
$v_t(W^{(i)})\subseteq W^{(i+rp)}$ for $v\in V^{(p)}$ for some $r\in \Z$. \\

\noindent
{\bf 2.3.3  Moose-Seiberg and Huang's Verlinde formula}\\
By the assumption, $W^i\boxtimes W^j$ decomposes into the direct 
sum of simple modules:
$$W^i\boxtimes W^j\cong \oplus_k (\underbrace{W^k\oplus \cdots \oplus W^k}_{N_{i,j}^k}).$$ 
We call $N_{ij}^k$ a fusion rule. 
A mysterious property of $C_2$-cofinite VOA whose modules are all 
completely reducible is a relation between 
the fusion rules and the entries of the $S$-matrix of $V$, 
which was mentioned by Verlinde \cite{V} and proved by 
Moose-Seiberg \cite{MoS} and Huang \cite{H2}. 
For example, the following is Corollary 5.4 in \cite{H2}. \\

\noindent
{\bf Theorem} \qquad {\it Let $V$ be a $C_2$-cofinite VOA of CFT-type 
and assume that $V\cong V'$ and all modules are completely reducible. 
Then $S$ is symmetric and the square $S^2$ is a permutation matrix 
which shifts $i$ to $i'$, where 
$W^0=V$ and $W^{k'}=(W^k)'$. Moreover, we have: 
$$\mbox{}\qquad \qquad \qquad 
N_{i,j}^k=\sum_{h=0}^s \frac{S_{ih}S_{jh}S_{hk'}}{S_{0h}} 
\qquad \qquad \hfill \mbox{{\rm (Verlinde Formula)}}.$$}

\section{Proof of Theorem A} 
We will first show that $V^G$ is semi-rigid. 
From the Main theorem in \cite{Mi3}, it is 
enough to show that $S(\Psi_{V^G})$ is a linear combination of 
trace functions on $V^G$-modules. 
For a finite cyclic automorphism group $<g>$, we have already proved 
it by Theorem 4 in \cite{Mi3}. 
Let ${\rm Irr}(G)$ be the set of irreducible $G$-characters and 
$M_{\chi}$ denotes 
an irreducible $G$-module affording to $\chi\in {\rm Irr}(G)$. 
Viewing $V$ as a $G$-module, 
$$V=\oplus_{\chi\in {\rm Irr}(G)}M_{\chi}\otimes V^{\chi}$$
for some space $V^{\chi}$. Clearly, $V^{\chi}$ is a $V^G$-module. 
Let us consider trace function $\Psi_V(gv)$ of $gv$ on $V$ for $g\in G$ 
and $v\in V^G$. 
$$\Psi_V(g v)=\sum_{\chi} \Tr_{M_{\chi}\otimes V^{\chi}} g o(v)q^{L(0)-c/24}
=\sum_{\chi} \chi(g)\Psi_{V^{\chi}}(v).$$
Hence we have 
$$S(\Psi_V(vg))=\sum_{\chi} \chi(g)S(\Psi_{V^{\chi}}(v))$$
The left side is a linear combination of trace functions on $V^g$-modules, 
which is clearly a linear combination of trace functions on $V^G$-modules. 
The sum of the left side for all $g\in G$ is equal to 
$|G|S(\Psi_{V^G}(v))$. 
Therefore, $S(\Psi_{V^G}(v)$ is a linear combination of trace functions 
on $V^G$-modules.

We next show that $V^G$-modules are completely reducible 
for a solvable group $G$.  
For a cyclic group $<g>$, $V^g$ is rational by 
Proposition 5 in \cite{Mi3}. Since $C_G(g)$ 
acts on $V^g$, we are able to obtain that 
$V^A$ is rational for an abelian normal subgroup $A$ of $G$. 
Repeating these steps, we have that $V^G$ is rational. \\
This completes the proof of Theorem A.

\section{Proof of Theorem B}
We will show that $V_L^{\sigma}$ is $C_2$-cofinite for any 
triality automorphism $\sigma$ of $L$. 
Abusing the notation, we use the same notation to denote 
an automorphism of $V_L$ lifted from $\sigma$.  

\subsection{Preliminary results}
Let's explain the definition of lattice VOA $V_L$, 
but we will give only necessary properties for the proof of Theorem B. 
See \cite{FLMe} for the precise definition. 
We first reduce the case. 
One advantage of proving $C_2$-cofiniteness is unexacting. For example, 
it is enough to show that $V_H^{\sigma}$ is $C_2$-cofinite 
for some $\sigma$-invariant full-sublattice $H$ of $L$, 
since $V_L^{\sigma}$ is a $V_H^{\sigma}$-module with 
a composition series of finite length, see Proposition \ref{prn:C2}. 

Since every lattice VOA is $C_2$-cofinite and $V^1\otimes V^2$ 
is $C_2$-cofinite when the both $V^i$ are $C_2$-cofinite, 
it is enough to prove the $C_2$-cofiniteness of $V_L^{\sigma}$ for 
$L=\Z x+\Z y$ with $y=\sigma(x)$ and $-x-y=\sigma^2(x)$. 
By taking a sublattice, we may assume that 
$\langle x,x\rangle=-2\langle x,y\rangle=18M>72$ and $M$ is even.   

Let's show a construction of a lattice VOA $V_L$ for $L=\Z x+\Z y$. 
Using a $2$-dimensional vector space $\C L=\C\otimes_{\Z}L$ with 
a nonsingular 
bilinear form $\langle\cdot,\cdot\rangle$, we first define 
a VOA $M_2(1):=S(\C L\otimes \C[t^{-1}]t^{-1})$ of free bosonic Fock space 
as follows:  
Consider $\C L$ as a commutative Lie algebra 
with a non-degenerated symmetric bilinear form, we construct 
the corresponding affine Lie algebra $\C L[t,t^{-1}]\oplus \C$ with a product 
$$  [v\otimes t^n, u\otimes t^m]=\delta_{n+m,0}n\langle v,u\rangle, $$
for $v,u\in \C L$.  Hereafter we use $v(n)$ to denote $v\otimes t^n$.  
We then consider its universal enveloping algebra  $U(\C L[t,t^{-1}])$. 
It has a commutative subalgebra $U(\C L[t])$.  
Using a one dimensional $U(\C L[t])$-module $\C e^{\gamma}$ with 
$$\mu(0)e^{\gamma}=\langle \mu,\gamma\rangle e^{\gamma}, \quad 
\mbox{ and }\quad \mu(n)e^{\gamma}=0 \mbox{ for }n>0  $$
for $\mu\in \C L$, we define a 
$U(\C L[t,t^{-1}])$-module: 
$$M_2(1)e^{\gamma}:=U(\C L[t,t^{-1}])\otimes_{U(\C L[t])}\C e^{\gamma}.$$ 
As vector spaces, $M_2(1)e^{\gamma}\cong S(\C L[t^{-1}]t^{-1})$ a 
symmetric tensor algebra. 
Since $L=\Z x+\Z y$, $M_2(1)e^{\gamma}$ is spanned by 
$$\{x(-i_1)\cdots x(-i_k)y(-j_1)\cdots y(-j_h)e^{\gamma}
\mid i_1\geq \cdots\geq i_k>0, j_1\geq\cdots\geq j_h>0\} $$
and we define its weight by $\sum i_s+\sum j_t+\frac{\langle \gamma,\gamma\rangle}{2}$. 

Forgetting about the vertex operators, 
a lattice VOA associated with $L$ is defined as a vector space  
$$V_L=\oplus_{\gamma\in L} M_2(1)e^{\gamma}.$$
We introduce an $\N$-gradation on $V_L=\oplus_{m\in \N} (V_L)_m$ by weights 
and let $(V_L)_m$ be the space of elements with weight $m$.  
We next define a vertex operator $Y(u,z)=\sum u_mz^{-m-1}\in {\rm End}(V_L)[[z,z^{-1}]]$ 
of $u\in V_L$ as follows: \\
First, a vertex operator of $v(-1)e^0$ for $v\in \C L$ and that of $e^{\gamma}$ are defined by 
$$ \begin{array}{rl}
Y(v(-1)e^0,z):&=\sum_{n\in \Z} v(n)z^{-n-1} \qquad \qquad \mbox{ and}\cr
Y(e^{\gamma},z):&=E^-(-\gamma,z)E^+(-\gamma,z)e^{\gamma}z^{\gamma},
\end{array}$$
where 
$$\begin{array}{l}
 E^{\pm}(\gamma,z)=\sum_{n=0}^{\infty}\frac{1}{n!}(\sum_{n\in {\Z}_+}
\frac{\gamma(\pm n)}{\pm n}z^{\mp n})^n 
\in {\rm End}(V_L)[[z^{\mp 1}]],  \cr
e^{\gamma}e^{\mu}=e^{\gamma+\mu} \qquad \mbox{ and }\qquad z^{\gamma}e^{\mu}
=z^{\langle \gamma,\mu\rangle}e^{\mu}.
\end{array}$$
The vertex operators of other elements are inductively defined by 
using normal products 
$$(v(-m)u)_n=
\displaystyle{\sum_{i=0}^{\infty}}(-1)^i\binom{-m}{i}\{v(-m-i)
u_{n+i}\!-\!(-1)^mu_{-m+n-i}v(i)\} \eqno{(2)}$$
for $v\in \C L$ and $u\in M_n(1)e^{\gamma}$, where$\binom{-m}{i}
=\frac{(-m)(-m-1)\cdots(-m-i+1)}{i!}$, 
and then we extend them linearly.  
We will frequently use this normal product development $(2)$.  

From now on, we use notation $\1$ to denote $e^0$ and call it a Vacuum. 
Important properties of the Vacuum $\1$ are $v_n\1$=0 and 
$v_{-1}\1=v$ for any $n\geq 0$ and $v\in V$. 
We also denote $M_n(1)e^0$ by $M_n(1)$. 
Let $\{x^{\ast},y^{\ast}\}$ be a dual basis of $\C L$ for $\{x,y\}$. 
Then we have a Virasoro element 
$\omega=\frac{1}{2}\{x(-1)x^{\ast}(-1)\1+y(-1)y^{\ast}(-1)\1\}$. We denote the operator 
$\omega_m$ by $L(m-1)$. 
Important properties of a Virasoro element are 
$$(L(-1)v)_m=-mv_{m-1} \quad \mbox{ and }\quad L(0)v=\wt(v)v.$$

For a VOA $V$ and its module $W$, we set:
$$\begin{array}{l}
C_m(V):=<v_{-m}u \mid v,u\in \oplus_{m\geq 1}V_m >_{\C}\qquad \mbox{ and} \cr 
C_2(V)W:=<v_{-2}w \mid v\in V, w\in W>_{\C}.
\end{array}$$

\subsection{In the free bosonic Fock space}
Viewing $\C L$ as a $\C[\sigma]$-module, we have 
$\C L=\C a\oplus \C a'$ with 
$\sigma(a)=e^{2\pi \sqrt{-1}/3}a$ and $\sigma(a')=e^{-2\pi \sqrt{-1}/3}a'$ and $\langle a,a'\rangle=1$. 
Then $u=a(-i_1)\cdots a(-i_h)a'(-j_1)\cdots a'(-j_k)\1$ is $\sigma$-invariant if and only if 
$h-k\equiv 0 \pmod{3}$. We note that $\omega=a(-1)a'=a'(-1)a$ is the Virasoro element.    
Set 
$$\CP_k=U(\C L[t,t^{-1}])\C L[t]t^k,$$ 
which are left ideals of $U(\C L[t,t^{-1}])$. Then 
$\CP_0e^0=0$ and $\CP_1e^{\gamma}=0$ for $\gamma\in L$.
From now on, $\equiv$ denotes a congruence relation modulo $C_2(M_2(1)^{\sigma})M_2(1)$. 
First of all, we will note the following lemma which comes from the normal product (2). \vspace{-4mm}\\

\begin{lmm}\label{lmm:reduction} Let $k=1$ or $2$ and $v\in M_2(1)$. 
For any $n,m>0$, there are $\lambda_{s,t}\in \Q$ such that 
$$  a(-n)a'(-m)v-\sum_{s,t>0} \lambda_{s,t}(a(-s)a'(-t)\1)_{-k}v \in \CP_{1-k}v. $$
For any $l,m,n\in \Z$ with $n\geq k$, there are $\lambda_{s,t,u}\in \Q$ such that 
$$a(-l)a(-m)a(-n)v -\sum \lambda_{s,t,u}(a(-s)a(-t)a(-u)\1)_{-k}v \in \CP_{1-k}.$$ 
\end{lmm}

\pr 
We will prove only the second case. Since 
$$\begin{array}{rl}
(a(-i)\gamma)_{-k}v=&
\sum_{j=0}^{\infty}\binom{-i}{j}(-1)^{j}\{a(-i-j)\gamma_{-k+i}-(-1)^i\gamma_{-i-k-j}a(j)\}v\cr
\in & \sum_{j=0}^{\infty}\binom{-i}{j}(-1)^{j}a(-i-j)\gamma_{-k+i}v+\CP_{0}v,
\end{array}$$
there are $\lambda_i, \lambda_{ij}\in \Q$ such that 
$$\begin{array}{l}
(a(-s)a(-t)a(-u)\1)_{-k}\in \sum_{i=0}^{\infty}\lambda_i a(-s-i)\{a(-t)a(-u)\1\}_{-k+i}+\CP_{0} \cr
\mbox{}\qquad =\sum_{i=0}^{\infty}\lambda_{i,j} a(-s-i)\sum_{j=0}^{\infty}a(-t-j)\{a(-u)\1\}_{-k+i+j}+\CP_{0} \cr
\mbox{}\qquad =\sum_{i=0}^{\infty}\lambda_{i,j} a(-s-i)\sum_{j=0}^{\infty}a(-t-j)a(-u-k+1+i+j)+\CP_{0}.
\end{array}$$
We hence have 
$a(-l)a(-m)a(-n)v-(a(-l)a(-m)a(-n+k-1)\1)_{-k}v$ is equivalent to a $\Q$-linear combination of 
$$\{a(-s)a(-t)a(-u)v \quad \mbox{ with } 0<u<n,  s\geq l, t\geq m\}$$
modulo $\CP_{-k+1}v$. Iterating these steps, we can reduce them to the case where $u<k$.  
\prend

We next explain an expression, which we will use. 
As a spanning set of  $V_L^{\sigma}$, we usually use elements of the form
$$  \mu=\sum_{t=0}^2\sigma^t(a(-i_1)\cdots a(-i_h)a'(-j_1)\cdots a'(-j_k)e^{\gamma}) 
\qquad \mbox{ with } i_s,j_t>0, $$
but we will permit to use $a(0)$ and $a'(0)$ so that 
$$ \mu=\frac{1}{\langle a,\gamma\rangle^s \langle b,\gamma\rangle^t}a(-i_1)\cdots a(-i_h)a'(-j_1)
\cdots a'(-j_k)a(0)^sa'(0)^t(\sum_{i=0}^2e^{\sigma^i(\gamma)}), \eqno{(3)}$$ 
where $ a(-i_1)\cdots a(-i_h)a'(-j_1)\cdots a'(-j_k)a(0)^sa'(0)^t$ 
is $\langle \sigma\rangle$-invariant and $a'(m)^n$ denotes 
$\underbrace{a'(m)\cdots a'(m)}_n$. 
From now on, $E^{\gamma}$ denotes $\sum_{i=0}^2e^{\sigma^i(\gamma)}$ 
and we will call $h$ and $k$ the numbers of $a$-terms and $a'$-terms, respectively.

\subsection{Modulo $C_2(M_2(1)^{\sigma})$}
For $u=a(-i_1)\cdots a(-i_h)a'(-j_1)\cdots a'(-j_k)\1$, if at least one of 
$\{i_s, j_t\mid s=1,...,h, t=1,...,k\}$ is not $1$, then we call $u$ a weight loss element.
Set 
$$\CS_2=\{a(-1)^{i}a'(-1)^j\1, a(-i)a(-j)a, a'(-i)a'(-j)a', a(-i)a'(-1)^2a, a(-i)a', \1 \mid i,j\in \N\}.$$ 

\begin{prn}\label{prn:reduction} 
Let $u=a(-i_1)\cdots a(-i_h)a'(-j_1)\cdots a'(-j_k)\1$ be a weight loss element.\\ 
If $|h-k|\geq 4$, then $u \in C_2(M_2(1)^{\sigma})$. \\
If $h-k=3$, then $u\in <a(-i_1)a(-i_2)a\mid i_1,i_2\in \N>_{\C}+\C_2(M_2(1)^{\sigma})$. \\  
If $h=k$, then $u\in <a(-\wt(u)+3)a'(-1)^2a, a(-\wt(u)+1)a'>_{\C}+\C_2(M_2(1)^{\sigma})$. \\
In particular, we have \quad$\displaystyle{M_2(1)^{\sigma}=C_2(M_2(1)^{\sigma})+<\CS_2>_{\C}}$.
\end{prn}

\pr 
We will prove the last statement. The others come from the same arguments. 
We first note that $\omega_0\beta=\beta_{-2}\1\in C_2(M_2(1)^{\sigma})$ for 
$\beta\in M_2(1)^{\sigma}$ and $a(-h)a(-k)a(-m)\1$ is congruent to a linear sum of 
elements of type $a(-i)a(-j)a$ since 
$\omega_0(a(-r_1)\cdots a(-r_k)\1)=\sum_{i=1}^k r_i a(-r_1)\cdots a(-r_i-1)\cdots a(-r_k)\1$.  
Suppose $h-k\geq 4$ and $u\not\in C_2(M_2(1)^{\sigma})+<\CS_2>_{\C}$. 
We take $u$ such that the total number $h+k$ is minimal. 
At least one of $i_s,j_t$ is not $1$. Since $h\geq 4$, by using suitable triple terms of $a$, 
we may assume $i_1=1$ by Lemma \ref{lmm:reduction}. Then by choosing other suitable triple $a$-terms, 
we may also assume $i_2=1$ and then $i_3=2$. 
Then 
$$2u -(a(-1)a(-1)a)_{-2}a(-i_4)\cdots a(-i_h)a'(-i_1)\cdots a'(-i_k)\1$$
is congruent to a linear sum of elements with the total number of terms is less than $h+k$, 
which contradicts the choice of $u$. 
We next treat the case $h-k=3$. By applying the same arguments to $a(-n)a'(-m)$, 
we can reduce to the case $h=3$ and $k=0$ as we desired.  
If $h=k$ and $h\geq 3$, then using the same argument as above, we can reduce to 
$u=a(-n)a'(-m)a(-1)a'(-1)$ and $n\geq 2$. 
If $m\geq 2$, then $u$ is congruent to a linear sum of $a(-n-m+1)a'(-1)^2a$ and $a(-n-m-1)a'$  
by Lemma \ref{lmm:reduction}. 
Therefore we obtain $M_2(1)^{\sigma}=C_2(M_2(1)^{\sigma})+<\CS_2>_{\C}$.
\prend

\subsection{A subring}
We note that 
$M_2(1)^{\sigma}/C_2(M_2(1)^{\sigma})$ is a commutative ring 
with the $-1$-normal product since $[\alpha_{-1},\beta_{-1}]
=\sum_{i=0}^{\infty}(-1)^i (\alpha_i\beta)_{-2-i}$ 
for any $\alpha,\beta\in V_L$.  
Let $\CO$ be the subspace of $M_2(1)^{\sigma}/C_2(M_2(1)^{\sigma})$ spanned by 
elements with the same number of $a$-terms and $b$-terms and 
$\CO^{even}$ the subspace of $\CO$ spanned by elements with even weights. 
Clearly, $\CO$ and $\CO^{even}$ are subrings of $M_2(1)^{\sigma}/C_2(M_2(1)^{\sigma})$ 
since $\wt(\alpha_{-1}\beta)=\wt(\alpha)+\wt(\beta)$.  
Let's study an algebraic structure of $\CO^{even}$. 

Set $\gamma(n)=a(-n+1)a'$. To simplify the notation, we sometimes omit 
subscript $-1$ denoting $-1$-normal product, for example, 
$\gamma(n)\gamma(m)$ denotes $\gamma(n)_{-1}\gamma(m)$. From 
$0\equiv \omega_0(a(-n)a'(-m)\1)=na(-n-1)a'(-m)\1+ma(-n)a'(-m-1)\1$, we have: \vspace{-4mm}\\

\begin{lmm}\label{lmm:aa}
$\mbox{}\qquad\displaystyle{ a(-n)a'(-m-1)\1\equiv \binom{-n}{m}\gamma(n+m+1)} \pmod{\omega_0V_L}\hfill {\rm (4)}$
\end{lmm}

\begin{prn}\label{prn:aaaa} 
$$\begin{array}{l}
a(-r)a(-m)a'(-n)a'\equiv \binom{-m}{n-1}\gamma(r+1)\gamma(m+n)
-\frac{(-1)^{n-1}(r+m+n-1)!(m+n+r+1)}{(r-1)!(m-1)!(n-1)!(m+1)(r+n)}\gamma(t)
\end{array}$$
modulo $\omega_0V_L$, 
where $t=r+m+n+1$. In particular, by replacing $r$ with $m$, we have 
$$\gamma(n+3)\equiv \frac{6}{(n-1)(n-2)(n+3)}\{\gamma(3)_{-1}\gamma(n)
-(n-1)\gamma(2)_{-1}\gamma(n+1)\}$$
for $n\geq 3$ and so $\gamma(n)\in C_1(M_2(1)^{\sigma})$ for $n\geq 6$. 
\end{prn}

\pr
The assertion comes from the direct calculation: 
$$\begin{array}{l}
\binom{-m}{n-1}\gamma(r+1)\gamma(m+n)\equiv (a(-r)a')_{-1}a(-m)a'(-n)\1\cr
\equiv \sum_i\binom{-r}{i}(-1)^i\{a(-r-i)a'(-1+i)-(-1)^{-r}a'(-r-1-i)a(i)\}
a(-m)a'(-n)\1 \cr
\equiv a(-r)a'(-1)a(-m)a'(-n)\1+\binom{r+m}{m+1}ma(-r-m-1)a'(-n)\1 \cr
\mbox{}\qquad -(-1)^r\binom{r+n-1}{n}a'(-r-1-n)na(-m)\1 \cr
\equiv a(-r)a(-m)a'(-n)a'+\{\binom{r+m}{m+1}m\binom{-r-m-1}{n-1}-(-1)^{n}\binom{r+n-1}{n}n
\binom{m+r+n-1}{r+n}\}\gamma(t) \cr
\equiv a(-r)a(-m)a'(-n)a'+
\frac{(-1)^{n-1} (r+m+n-1)!(m+r+n+1)}{(r-1)!(m-1)!(n-1)!(m+1)(r+n)}\gamma(t). \hfill \mbox{\prend}
\end{array}$$

For example, we will use the following: 
$$\begin{array}{c}
2\gamma(6)\equiv \gamma(3)\gamma(3)-2\gamma(2)\gamma(4), \qquad 
7\gamma(7)\equiv \gamma(3)\gamma(4)-3\gamma(2)\gamma(5), \cr
16\gamma(8)\equiv \gamma(3)\gamma(5)-4\gamma(2)\gamma(6), \qquad
30\gamma(8)\equiv \gamma(4)\gamma(4)-6\gamma(2)\gamma(6). 
\end{array}$$
\vspace{-2mm}

\begin{lmm}\label{lmm:CO}\qquad
$\CO=<\gamma(2)^n, \gamma(n+1), \gamma(2)\gamma(m)\1\mid n,m=2,\ldots >_{\C}$.
\end{lmm}

\pr 
By Proposition \ref{prn:reduction}, $\CO$ is spanned by $\{a(-1)^na'(-1)^n\1, a(-n)a'(-1)^2a, \gamma(m)\}$.
By Proposition \ref{prn:aaaa}, we get $a(-n)a'(-1)^2a-\gamma(2)\gamma(n+1)\in \Q\gamma(n+3)$.  
We also have that $a(-1)^na'(-1)^n\1-\gamma(2)^n$ is congruence to a linear 
sum of $a(-2n+3)a'(-1)^2a$ and $\gamma(2n)$ modulo $C_2(M_2(1)^{\sigma})$, which 
proves the desired result. 
\prend

Set $\CS_1=\{a(-i_1)a(-i_2)a, a'(-i_1)a'(-i_2)a', a(-i_3)a', \1 \mid i_1,i_2\leq 5, i_3\leq 4 \}$. 
\vspace{-4mm}\\

\begin{prn}\label{prn:C1}  
$M_2(1)^{\sigma}=C_1(M_2(1)^{\sigma})+<\CS_1>_{\C}$. 
In particular, $M_2(1)^{\sigma}$ is $C_1$-cofinite. 
\end{prn}

\pr 
To simplify the notation, set $C_1=C_1(M_2(1))^{\sigma}$ in this proof. 
Suppose that the proposition is false and let 
$$u=a(-i_1)\cdots a(-i_h)a'(-j_1)\cdots a'(-j_k)\1\not\in C_1+<\CS_1>_{\C}.$$ 
We take $u$ such that the number of terms is minimal. 
By Lemma \ref{lmm:reduction} and \ref{lmm:aa}, we may assume 
$u=a(-i_1)a(-i_2)a$  or $u=a(-m)a'$. 
By Lemma \ref{lmm:aa} and Proposition \ref{prn:aaaa}, we obtain $a(-m)a'\in C_1$ for $m\geq 5$. 
Since $C_1$ is closed by the $0$-th product, we have:
$$\begin{array}{l}
(1) \qquad C_1\ni (a(-k+1)a')_0(a(-1)a(-1)a)=3(k-1)a(-k)a(-1)^2\1  \qquad \mbox{ and so}\cr
(2) \qquad C_1\ni (a(-n)a')_{0}a(-1)^2a(-k)\1=2a(-n-1)a(-k)a+ka(-n-k)a(-1)a   
\end{array}$$
for $k\geq 6$ and any $n$. \prend

We next express $\CO$ as a $\C[\gamma(2)]$-module. We need the following lemma. \vspace{-4mm}\\

\begin{lmm}\label{lmm:gamma8}
$\mbox{}\qquad 120\gamma(7)\1\equiv 8\gamma(2)\gamma(5)\1+\gamma(2)^2\gamma(3)\1$ \\
$\mbox{}\qquad \qquad \qquad \qquad 60\gamma(8)\1\equiv 6\gamma(2)\gamma(3)^2\1-13\gamma(2)^2\gamma(4)\1$.
\end{lmm}

\pr 
Since $\mbox{}\quad 
0\equiv (a(-1)a(-1)a)_{-2}a'(-1)a'(-1)a' $\\
$\equiv 3a(-1)^2a(-2)a'(-1)^2a'+18a(-4)a(-1)a'(-1)a'+18a(-3)a(-2)a'(-1)a'+18\gamma(7)$, \\
we have:
$$a(-1)^2a(-2)a'(-1)^2a'\equiv -6a(-4)a(-1)a'(-1)a'-6a(-3)a(-2)a'(-1)a'-6\gamma(7).$$
Using Proposition \ref{prn:aaaa} and the above lemma, we obtain the first congruence expression:
$$\begin{array}{rl}
\multicolumn{2}{l}{\gamma(2)^2\gamma(3)\equiv(a(-1)a')_{-1}\{a(-1)a'(-1)a(-2)a'\}
+\gamma(2)\{2a(-4)a'+a'(-3)a(-2)\1\}}\cr
\equiv& (a(-1)a'(-1)a(-1)a'(-1)a(-2)a'+a(-3)a'(-1)a(-2)a'\cr
&+2a(-4)a(-1)a'(-1)a'+2a'(-3)a(-1)a(-2)a'+5\gamma(2)\gamma(5) \cr
\equiv& -6a(-4)a(-1)a'(-1)a'-6a(-3)a(-2)a'(-1)a'-6\gamma(7) \cr
&+a(-3)a'(-1)a(-2)a'+2a(-4)a'(-1)^2a+2a'(-3)a(-1)a(-2)a'+5\gamma(2)\gamma(5)\cr
\equiv& -4a(-4)a(-1)a'(-1)a'-5a(-3)a(-2)a'(-1)a'-6\gamma(7)\cr
&+2a'(-3)a(-1)a(-2)a'+5\gamma(2)\gamma(5) \cr
\equiv& -4\{\binom{-4}{0}\gamma(2)\gamma(5)-[4+\binom{-4}{2}]\gamma(7)\}
-5\{3\gamma(2)\gamma(5)-28\gamma(7)\}-6\gamma(7) \cr
&+2\{\binom{-2}{2}\gamma(2)\gamma(5)-[2\binom{-4}{2}+3\binom{-2}{4}]\gamma(7)\}
+5\gamma(2)\gamma(5) \cr
\equiv& 120\gamma(7)-8\gamma(2)\gamma(5). 
\end{array}$$
By expanding $0\equiv (a(-1)a(-1)a)_{-2}a'(-2)a'(-1)a'$, we have, 
$$\begin{array}{rl}
\multicolumn{2}{l}{-(a(-2)a(-1)a(-1))a'(-2)a'(-1)a' }\cr
\equiv & 4a(-4)a(-1)a'(-2)a' +4a(-3)a(-2)a'(-2)a' +4a(-5)a(-1)a'(-1)a' \cr
 &+4a(-4)a(-2)a'(-1)a' +2a(-3)a(-3)a'(-1)a'+8\gamma(8)
\end{array}$$
and then we obtain: 
$$\begin{array}{l}
2\gamma(2)\gamma(2)\gamma(4)\equiv -(a(-1)a')_{-1}\{a(-1)a(-2)a'(-2)a'-16\gamma(6))\} \cr
\equiv -(a(-1)^2a'(-1)^2ba(-2)a'(-2)-a(-3)a(-2)a'(-2)a'-2a(-4)a(-1)a'(-1)a'(-2) \cr
\mbox{}\quad -a'(-3)a(-1)a(-2)a'(-2)-2a'(-4)a(-1)a'(-1)a(-2)+16\gamma(2)\gamma(6) \cr
\equiv 4a(-4)a(-1)a'(-2)a' +4a(-3)a(-2)a'(-2)a' +4a(-5)a(-1)a'(-1)a' \cr
\mbox{}\quad+4a(-4)a(-2)a'(-1)a'+2a(-3)a(-3)a'(-1)a'+8\gamma(8)\cr
\mbox{}\quad-a(-3)a'(-1)a(-2)a'(-2)-2a(-4)a(-1)a'(-1)a'(-2)-a'(-3)a(-1)a(-2)a'(-2)\cr
\mbox{}\quad-2a'(-4)a(-1)a'(-1)a(-2)+16\gamma(2)\gamma(6) \cr
\equiv 2a(-3)a(-2)a'(-2)a' +4a(-5)a(-1)a'(-1)a' +4a(-4)a(-2)a'(-1)a'   \cr
\mbox{}\quad+2a(-3)a(-3)a'(-1)a'+16\gamma(2)\gamma(6)+8\gamma(8) \cr
\equiv 2(180\gamma(8)-3\gamma(3)\gamma(5))+4(\gamma(2)\gamma(6)-20\gamma(8)) \cr
\mbox{}\quad+4(\gamma(3)\gamma(5)-64\gamma(8))+2(-90\gamma(8)+\gamma(4)\gamma(4))
+16\gamma(2)\gamma(6)+8\gamma(8) \cr
\equiv -120\gamma(8)+12\gamma(2)\gamma(3)^2-24\gamma(2)^2\gamma(4). \prend
\end{array}$$
By the above lemma, the direct calculation shows:
$$\begin{array}{l}
2\gamma(4)\gamma(4)\equiv 12\gamma(2)\gamma(6)+60\gamma(8)
\equiv 12\gamma(2)\gamma(3)^2-25\gamma(2)^2\gamma(4), \cr
15\gamma(3)\gamma(5)=60\gamma(2)\gamma(6)+240\gamma(8)
\equiv 54\gamma(2)\gamma(3)^2-112\gamma(2)^2\gamma(4), \cr
120\gamma(3)\gamma(4)\equiv 120(7\gamma(7)+3\gamma(2)\gamma(5))
\equiv 7\gamma(2)^2\gamma(3)+416\gamma(2)\gamma(5). 
\end{array}$$
Therefore, $\CO^{even}$ has a subring 
$\displaystyle{\CO_{\Q}^{even}=\Q[\gamma(2)]\gamma(2)+\Q[\gamma(2)]\gamma(3)\gamma(3)
+\Q[\gamma(2)]\gamma(4)}$. 

\subsection{Elements $a(-1)a(-1)a$}
We denote $a(-1)a(-1)a$ and $a'(-1)a'(-1)a'$ by $\alpha$ and $\beta$, respectively. \vspace{-4mm}\\

\begin{lmm}\label{lmm:gamma3}
$\mbox{}\qquad \gamma(2)_{-1}\gamma(2)_{-1}\gamma(2)
\equiv \alpha_{-1}\beta-264\gamma(2)_{-1}\gamma(4)\1+117\gamma(3)_{-1}\gamma(3)$
\end{lmm}

\pr  From the direct calculation, we have:
$$\begin{array}{l}
\alpha_{-1}\beta\equiv (a(-1)a(-1)a)_{-1}a'(-1)a'(-1)a' \cr
\equiv a(-1)^3a'(-1)^3\1+18a(-3)a(-1)a'(-1)a'+9a(-2)a(-2)a'(-1)a'+18a(-5)a'.
\end{array}$$
Therefore, by Proposition \ref{prn:aaaa}, we obtain:
$$\begin{array}{rl}
\gamma(2)^3\equiv& (a(-1)a')_{-1}\{a(-1)a'(-1)a(-1)a'+2\gamma(4)\} \cr
\equiv& (a(-1)^3a'(-1)^3\1+2a(-3)a(-1)a'(-1)a'+2a'(-3)a(-1)a'(-1)a+2\gamma(2)\gamma(4) \cr
\equiv& \alpha_{-1}\beta-14\{\gamma(2)\gamma(4)-9\gamma(6)\}
-9\{\gamma(3)^2-16\gamma(6)\}-18\gamma(6)+2\gamma(2)\gamma(4) \cr
\equiv& \alpha_{-1}\beta-264\gamma(2)\gamma(4)+117\gamma(3)^2. \prend
\end{array}$$

\subsection{The action of $\gamma(4)$}
In this subsection, we will consider elements modulo $C_2(M_2(1)^{\sigma})$ and we abuse 
$=$ to denote $\equiv$. 
By \S 4.4, we have shown that $\CO_{\Q}^{even}$ is closed by the $-1$-product and 
$$\CA_{\Q}^{even}=\Q[\gamma(2)]\gamma(4)+\Q[\gamma(2)]\gamma(3)\gamma(3)$$
is an ideal modulo $C_2(M_2(1)^{\sigma}$. Let $\CQ$ be an ideal generated by 
$\alpha_{-1}\beta$. We note $\alpha_{-1}\beta\equiv \gamma(2)^3+264\gamma(2)\gamma(4)-117\gamma(3)^2$. 
We will see the action of $\gamma(4)$ on $\CA_{\Q}^{even}$. \vspace{-4mm}\\

\begin{lmm}\label{lmm:CQ}
$\CQ=\CO_{\Q}^{even}$.
\end{lmm}

\pr 
We already know $\gamma(4)^2\equiv 6\gamma(2)\gamma(3)^2-\frac{25}{2}\gamma(2)^2\gamma(4)$. 
Since $\gamma(3)\gamma(5)\equiv 54\gamma(2)\gamma(3)^2-112\gamma(2)^2\gamma(4)$, we have:
$$\begin{array}{rl}
1800\gamma(4)\gamma(3)^2\equiv &15\gamma(3)\{7\gamma(2)^2\gamma(3)+416\gamma(2)\gamma(5)\} \cr
\equiv &105\gamma(2)^2\gamma(3)^2+416\gamma(2)\{54\gamma(2)\gamma(3)^2-112\gamma(2)^2\gamma(4)\}\cr
\equiv &22569\gamma(2)^2\gamma(3)^2-46592\gamma(2)^3\gamma(4).
\end{array}$$
Therefore the action of $\gamma(4)$ on $\CA_{\Q}^{even}$ is expressed by 
$\displaystyle{\gamma(2)^2 \left( \begin{array}{cc} 
\frac{22569}{1800} & \frac{-46592}{1800} \cr
& \cr
6 & \frac{-25}{2} \end{array}\right)}$. 
The eigenpolynomial of $1800\gamma(4)$ is $X^2-69X-4608900$  
and its discriminant is $3\sqrt{2048929}$, which is not a rational number. 
Therefore, the action of $\gamma(4)/\gamma(2)^2$ 
on $\Q\gamma(2)\gamma(4)+\Q\gamma(3)^2$ is irreducible over $\Q$. 
Furthermore, since 
$$\begin{array}{l}
(\alpha_{-1}\beta)_{-1}\gamma(4)\equiv(\gamma(2)^3-264\gamma(2)\gamma(4)-117\gamma(3)^2)\gamma(4) \cr
\mbox{}\quad \equiv \gamma(2)^3\gamma(4)-264\gamma(2)\{6\gamma(2)\gamma(3)^2
-\frac{25}{2}\gamma(2)^2\gamma(4)\}\cr
\mbox{}\qquad -117\gamma(3)\{\frac{7}{120}\gamma(2)^2\gamma(3)+\frac{416}{120}\gamma(2)\gamma(5)\} \cr
\mbox{}\quad \equiv 3301\gamma(2)^3\gamma(4)-\{1584+\frac{273}{40}\}\gamma(2)\gamma(3)^2
-\frac{39\times 52}{5}\gamma(2)\{\frac{54}{15}\gamma(2)\gamma(3)^2-\frac{112}{15}\gamma(2)^2\gamma(4)\} \cr
\mbox{}\quad \equiv (3301+\frac{13\times 52\times 112}{25})\gamma(2)^2\gamma(4)-\{1584+\frac{273}{40},
+\frac{39\times 52\times 18}{25}\}\gamma(2)^2\gamma(3)^2,  
\end{array}$$
we have $\CQ_{\Q}^{even}\cap \CA_{\Q}^{even}\not=0$ and so 
$$(<\alpha_{-1}\beta, \gamma(4)\alpha_{-1}\beta,\gamma(4)^2\alpha_{-1}\beta>_{\Q})_n=(\CO_{\Q}^{even})_n
\qquad \mbox{ for }n\geq 14. \hfill \mbox{\prend}$$

\subsection{Nilpotency of $\alpha$ modulo $C_2(V_L^{\sigma})V_L$}
From now on, $\equiv$ denotes the congruence modulo $C_2(V_L^{\sigma})V_L$.  We next show that \vspace{-4mm} \\

\begin{lmm}\label{lmm:nilpotent}
  $\mbox{}\qquad(a(-i_1)a(-i_2)a)_{-1}$ and $(a'(-j_1)a'(-j_2)a')_{-1}$ are all nilpotent in \\
$M_2(1)^{\sigma}/(C_2(V_L^{\sigma})\cap M_2(1))$ for any $i_1,i_2,j_1,j_2$. 
\end{lmm}

\pr
Except $\alpha$ and $\beta$, the square of the remainings are zero by Proposition \ref{prn:reduction}. 
We will prove that $\alpha_{-1}$ is nilpotent. 
Since $\wt(e^x)=9M$, $\wt(e^{x-y})=27M$ and $\wt(e^{2x+y})=27M$ 
for $y=\sigma(x)$, we have $e^y_{-1-k}e^{-x}=e^{-x-y}_{-1-k}e^{x}=0$ for $k<9M$ and so   
$$E^{x}_{-1-k}E^{-x}=\sum_{i=0}^2\sigma^i(E^x_{-1-k}e^{-x})
=\sum_{i=0}^2\sigma^i(e^x_{-1-k}e^{-x})\in M_2(1)^{\sigma}\cap C_2(V_L^{\sigma}) \quad 
\mbox{ for } 1<k<9M, $$
where $E^x$ denotes $e^x+e^y+e^{-x-y}$. 
Multiplying $(\alpha_{-1})^{6M+9}$ to $E^{x}_{-4}e^{-x}$, 
the number of $a$-terms in $(\alpha_{-1})^{6M+9}E^{x}_{-4}e^{-x}$ is 
$6$ more than that of $a'$-terms and so all elements with weight loss vanished. Hence  
$$(\alpha_{-1})^{6M+9}E^{x}_{-4}e^{-x}\equiv 
\frac{1}{(18M+3)!}(\alpha_{-1})^{6M+9}(x(-1))^{18M+3}\1\in C_2(V_L^{\sigma}).$$
Set $x=r a+s a'$, then since we multiply many $a(-1)$, 
$(\alpha_{-1})^{6M+9+k}$ annihilates all elements except for $a(-1)$ and $a'(-1)$ by 
Proposition \ref{prn:reduction} and so we have:
$$\begin{array}{l}
(\alpha_{-1})^{6M+9}(x(-1))^{18M+3}\1\equiv a(-1)^{18M+27}(r a(-1)+s a'(-1))^{18M+3}\1 \cr
\mbox{}\quad \qquad\equiv \sum_{i=0}^{18M+3}\binom{18M+3}{i}r^{18M+3-i}s^i a(-1)^{36M+30-i}\gamma(2)^i \cr
\mbox{}\quad\qquad\equiv 
\sum_{i=0}^{6M+1}\binom{18M+3}{3i}r^{18M+3-3i}s^{3i}(\alpha_{-1})^{12M+10-i}\gamma(2)^{3i}\cr
\mbox{}\qquad\qquad+\sum_{i=0}^{6M}\binom{18M+3}{3i+1}r^{18M+2-3i}s^{3i+1}(\alpha_{-1})^{12M+9-i}a(-1)a(-1)\gamma(2)^{3i+1}\cr
\mbox{}\qquad\qquad+\sum_{i=0}^{6M}\binom{18M+3}{3i+2}r^{18M+1-3i}s^{3i+2}(\alpha_{-1})^{12M+9-i}a(-1)\gamma(2)^{3i+2}.
\end{array}$$
Similarly, since we obtain  
$$\begin{array}{l}
(\alpha_{-1})^{6M+9}E^{x}_{-4}a(-1)e^{-x}=
\alpha_{-1}^{6M+9}(x(-1))^{18M+3}a(-1)\1+\alpha_{-1}^{6M+9}\langle a,x\rangle(x(-1))^{18M+4}\1 \cr
\mbox{}\qquad=\alpha_{-1}^{6M+9}(x(-1))^{18M+3}a(-1)\1+\alpha_{-1}^{6M+9}\langle a,x\rangle 
E^x_{-5}e^{-x} \cr
\mbox{}\qquad \equiv \alpha_{-1}^{6M+9}(x(-1))^{18M+3}a(-1)\1 \cr
\mbox{}\qquad \equiv
\sum_{i=0}^{6M+1}\binom{18M+3}{3i}r^{18M+3-3i}s^{3i}\alpha_{-1}^{12M+10-i}a(-1)\gamma(2)^{3i}\cr
\mbox{}\qquad\quad+\sum_{i=0}^{6M}\binom{18M+3}{3i+1}r^{18M+2-3i}s^{3i+1}\alpha_{-1}^{12M+9-i+1}\gamma(2)^{3i+1}\cr
\mbox{}\qquad\quad+\sum_{i=0}^{6M}\binom{18M+3}{3i+2}r^{18M+1-3i}s^{3i+2}\alpha_{-1}^{12M+9-i}a(-1)^2\gamma(2)^{3i+2}\cr
\multicolumn{1}{l}{\mbox{and}} \cr
\alpha_{-1}^{6M+9}E^{x}_{-4}a(-1)^2e^{-x}=\alpha_{-1}^{6M+9}(x(-1))^{18M+3}a(-1)^2\1
+2\langle a,x\rangle\alpha_{-1}^{6M+9}(x(-1))^{18M+4}a \cr
\mbox{}\qquad\quad+2\langle a,x\rangle^2\alpha_{-1}^{6M+9}(x(-1))^{18M+5}\1, \cr
\mbox{}\qquad\equiv \alpha_{-1}^{6M+9}(x(-1))^{18M+3}a(-1)^2\1 \cr
\mbox{}\qquad\equiv 
\sum_{i=0}^{6M+1}\binom{18M+3}{3i}r^{18M+3-3i}s^{3i}\alpha_{-1}^{12M+10-i}a(-1)^2\gamma(2)^{3i}\cr
\mbox{}\qquad\quad+\sum_{i=0}^{6M}\binom{18M+3}{3i+1}r^{18M+2-3i}s^{3i+1}\alpha_{-1}^{12M+9-i+1}a(-1)\gamma(2)^{3i+1}\cr
\mbox{}\qquad\quad+\sum_{i=0}^{6M}\binom{18M+3}{3i+2}r^{18M+1-3i}s^{3i+2}\alpha_{-1}^{12M+9-i+1}\gamma(2)^{3i+2},\cr
\end{array}$$
we have 
$$\begin{array}{l}
a(-1)\alpha_{-1}^{6M+9}(x(-1))^{18M+3}\1\in C_2(V_L^{\sigma})V_L,  \cr
a'(-1)\alpha_{-1}^{6M+9}(x(-1))^{18M+3}\1
\in C_2(V_L^{\sigma})V_L \qquad \mbox{ and}\cr
a(-1)a(-1)\alpha_{-1}^{6M+9}(x(-1))^{18M+3}\1\in C_2(V_L^{\sigma})V_L. 
\end{array}$$
Hence  
$$\alpha_{-1}^{6M+9+k}a(-1)^ea'(-1)^k(x(-1))^{18M+3}\1$$ 
is a linear sum of elements of the form 
$$\alpha_{-1}^{6M+9+k}v_{-1}(u\cdot(x(-1))^{18M+3}\1),$$ 
where $v$ is a $\sigma$-invariant element and 
$u \in \{ \1_{-1}, a(-1), a(-1)a(-1)\}$ 
by Lemma \ref{lmm:reduction}. Therefore we obtain  
$$\alpha_{-1}^{6M+9+k}a(-1)^ea'(-1)^k(x(-1))^{18M+3}\1\in C_2(V_L^{\sigma})V_L$$
for any $e, k\geq 0$. We also get a similar result for $y=\sigma(x)$ as for $x$. Therefore we have: 
$$\alpha_{-1}^{12M+18}(\lambda x(-1)+\mu y(-1))^{36M+6}\1\in C_2(V_L^{\sigma})V_L$$
for any $\lambda, \mu\in \C$. By choosing 
suitable $\lambda$ and $\mu$ so that $\lambda x(-1)+\mu y(-1)=a(-1)$, 
we have 
$$\alpha_{-1}^{12M+18}a(-1)^{36M+6}\1=\alpha_{-1}^{48M+24}\1 \in C_2(V_L^{\sigma}),$$
which implies that $\alpha_{-1}$ is nilpotent modulo $C_2(V_L^{\sigma})$. 
Similarly, $\beta_{-1}$ is nilpotent.
\prend

Since $\alpha, \beta$ are nilpotent and $\CO^{even}=\CO^{even}\alpha_{-1}\beta$, 
we have the following: \vspace{-4mm}\\

\begin{prn}\label{prn:M2C2}
$\mbox{}\qquad 
\displaystyle{\dim \left(M_2(1)^{\sigma}/(M_2(1)^{\sigma}\cap C_2(V_L^{\sigma}))\right) <\infty}$. 
\end{prn}

\subsection{$C_2$-cofiniteness of $V_L^{\sigma}$}

By the previous proposition, there is an integer $N_0$ such that 
$v^1_{-1}\cdots v^k_{-1}\gamma\in C_2(V_L^{\sigma})$ for any $v^i\in \CS_1$ and $\gamma\in V_L^{\sigma}$ 
if $\wt (v^1_{-1}\cdots v^k_{-1}\1) \geq N_0$. Set $N=N_0+9M+30$.  

Our final step is to prove that 
$$V_L^{\sigma}=C_2(V_L^{\sigma})+\oplus_{n\leq N}(M_2(1))_n^{\sigma} 
+\oplus_{n\leq N}(M_2(1)E^x)_n^{\sigma} 
+\oplus_{n\leq N}(M_2(1)E^{-x})_n^{\sigma},$$
which implies the $C_2$-cofiniteness of $V_L^{\sigma}$. For $\mu\not=0$, set \\
$\CR=\left\{d^k_{i_k}\cdots d^1_{i_1}b_{i_0}a(r)a'(0)E^{\mu}\mid 
\begin{tabular}{l}(i) $i_k\leq,\ldots,\leq i_1\leq -1$, $i_0\leq 0$, and \\
(ii) $d^i\in \CS_1$, $\wt(b_{i_0}a(r)a'(0)E^{\mu})-\wt(E^{\mu})\leq 30$
\end{tabular}\right\}$.

\begin{prn}\label{prn:spanningset}
$(M_2(1)E^{\mu})^{\sigma}=<\CR>_{\C}+C_2(V_L^{\sigma})$. 
In particular, if $v\in (M_2(1)E^{\mu})^{\sigma}$ has a weight 
greater than $\wt(E^{\mu})+N_0+30$, then $v\in C_2(V_L^{\sigma})$. 
\end{prn}

\pr
Suppose false and 
we take $u\not\in <\CR>_{\C}+C_2(V_L^{\sigma})$ such that $\wt(u)$ is minimal.  
Since $M_2(1)E^{\mu}$ is an irreducible $M_2(1)^{\sigma}$-module, 
we may assume 
$$ u=c^k_{i_k}\cdots c^1_{i_1}E^{\mu}$$
with $c^i\in M_2(1)^{\sigma}$. 
We take the above expression such that $\sum_{i=1}^k \wt(c^i)$ is minimal 
and if $\sum_{i=1}^k \wt(c^i)$ is the same, then $k$ is maximal.  
Since $(e_{-1}f)_k=\sum_{i=0}^{\infty}( e_{-1-i}f_{k+i}+f_{k-1-i}e_i)$ and $\wt(e_{-1}f)=\wt(e)+\wt(f)$, 
we may assume $c^i\in \CS_1$. Also, since $e_sf_t-f_te_s=\sum_{i=0}^{\infty}\binom{s}{i}(e_if)_{s+t-i}$ 
and $\wt(e_if)<\wt(e)+\wt(f)$ for $i\geq 0$, we may assume $i_k\leq \cdots \leq i_1$. 
By the minimality of $\wt(u)$, we have $0\leq i_k$ and 
$$\sum_{i=1}^k \wt(c^i)=(\wt(u)-\wt(E^{\mu}))+\sum_{i=1}^k (1+i_j).$$  
To simplify the notation, we will call 
$\sum_{i=1}^k(1+i_j)$ $\sigma$-loss weight. Since 
$\wt(E^{\mu})$ and $\wt(u)$ are fixed, we have chosen $u=c^k_{i_k}\cdots c^1_{i_1}E^{\mu}$ 
such that the $\sigma$-loss weight is minimal. 
We note that $\wt(c^i)\leq 11$ for $c^i\in \CS_1$.  
On the other hand, by Lemma \ref{lmm:reduction}, $u$ is also a linear sum of elements of the form 
$$e^r_{-1}\cdots e^1_{-1}F,$$
where $e^i\in M_2(1)^{\sigma}$ and $F$ is one of 
$$\CD=\{a(-m-n)a'(0)E^{\mu}, a(-m)a(-n)a(0)E^{\mu}, a'(-m)a'(-n)a'(0)E^{\mu}\}.$$ 
By the minimality of $\wt(u)$, $u$ is a linear sum of elements in $\CD$ and $m+n+\wt(E^{\mu})=\wt(u)$. 
We assert that the $\sigma$-loss weight of $u$ is less than or equal to three. 
For the elements $a(-m-n)a'(0)E^{\mu}$, 
we get $a(-m-n)a'(0)E^{\mu}=(a'(-m-n-1)a)_1E^{\mu}$, which has only 
$\sigma$-loss weight two. 
Before we start the proof for the remaining case, we note
$$\begin{array}{l}
(a'(-m-1)a)_1(a'(-n-1)a)_1E^{\mu}=(a'(-m-1)a)_1a(-n)a'(0)E^{\mu} \cr
\mbox{}\qquad =\sum\binom{-m-1}{i}(-1)^i(-1)^ma(-m-i)a'(i)a(-n)a'(0)E^{\mu} \cr
\mbox{}\qquad =\binom{-m-1}{n}(-1)^{n+m}na(-m-1-n)a'(0)E^{\mu}+(-1)^ma(-m)a'(0)a(-n)a'(0)E^{\mu}. \cr
\end{array}$$
Suppose $a(-m)a(-n)a(0)E^{\mu}$ has a $\sigma$-loss weight greater than three. 
By ignoring elements with $\sigma$-loss weight less than three, we have 
$$\begin{array}{l}
\frac{\langle b,\mu\rangle^2}{\langle a,\mu\rangle}a(-m)a(-n)a(0)E^{\mu}
=a(-m)a'(0)a(-n)a'(0)E^{\mu}\cr
\mbox{}\quad\equiv (a'(-m-\!1)a)_1(a'(-n\!-\!1)a)_1E^{\mu} \cr
\mbox{}\quad =(a'(-m-1)a)_1(a(-n-1)a')_1E^{\mu}+(a'(-m-1)a)_1(\omega_0\gamma(n+1)+\cdots))_1E^{\mu}\cr
\mbox{}\quad \equiv (a'(-m-1)a)_1(a(-n-1)a')_1E^{\mu}+(\omega_0\gamma(n+1))_1(a'(-m-1)a)_1E^{\mu}\cr
\mbox{}\quad \equiv (a'(-m-1)a)_1a'(-n)a(0)E^{\mu}-\gamma(n+1)_0(a'(-m-1)a)_1E^{\mu}\cr
\mbox{}\quad \equiv \sum \binom{-m-1}{i}(-1)^i\{a'(-m-1-i)a(1+i)-(-1)^{m+1}a(-m-i)a'(i)\}a'(-n)a(0)E^{\mu}\cr
\mbox{}\quad \equiv \binom{-m-1}{n-1}(-1)^{n-1}na'(-m-n)a(0)E^{\mu}-(-1)^{m+1}a(-m)a'(0)a'(-n)a(0)E^{\mu}\cr
\mbox{}\quad \equiv \lambda_1 a(-m-n)a'(0)E^{\mu}+\mu_1 a(-m)a'(-n)E^{\mu} \cr
\mbox{}\quad \equiv \lambda_2 a(-m-n)a'(0)E^{\mu}+\mu_2 a'(-m-n)a(0)E^{\mu}\equiv 0 \cr 
\end{array}$$
for some $\lambda_i$ and $\mu_j$, which is a contradiction. 
Therefore, the $\sigma$-loss weight of $u$ is less than or equal to three. 
In particular, $k\leq 3$ and $\wt(u)-E^{\mu}\leq 30$. Therefore, 
the elements $a(-m-n)a'(0)E^{\mu}$ and $\gamma(n+1)_0(a'(-m-1)a)_1E^{\mu}$ are all in $<\CR>_{\C}+C_2(V_L^{\sigma})$. 
In order to show $a(-m)a(-n)a(0)E^{\mu}\in <\CR>_{\C}+C_2(V_L^{\sigma})$, we have exactly the same 
congruence expressions as above modulo $<\CR>_{\C}+C_2(V_L^{\sigma})$. 
\prend

Set $K=M_2(1)^{\sigma}+(M_2(1)E^{x})^{\sigma}+(M_2(1)E^{-x})^{\sigma}$. 
Since we have already shown that if $v\in K$ and $\wt(v)>N$, then $v\in C_2(V_L^{\sigma})$. 
The remaining is to show 
$$ V_L^{\sigma}=K+C_2(V_L^{\sigma}).$$
By Proposition \ref{prn:spanningset}, it is enough to show that 
$$a(-n)a'(0)E^{\mu} \in K+C_2(V_L^{\sigma})$$
for $1\leq n \leq 30$ and $\mu\not\in \{0,\pm x,\pm y, \pm(x+y)\}$. We first 
treat the following case:\vspace{-4mm}\\

\begin{lmm}\label{lmm:modulo3}  For any $n+m\equiv 0 \pmod{3}$, we have $E^{mx+ny}\in C_2(V_L^{\sigma})+K$.
\end{lmm}

\pr  We note that if $n+m\equiv 0 \pmod{3}$, then there is $\gamma\in L$ such that 
$E^{mx+ny}=E^{\pm(\sigma(\gamma)-\gamma)}$. 
Set $2k=\langle \gamma,\gamma\rangle$. Then 
since $\langle \gamma-\sigma(\gamma),\gamma-\sigma(\gamma)\rangle=6k$, 
we have 
$$\begin{array}{l}
E^{\gamma}_{-1-k}E^{-\gamma}\in M_2(1)+E^{\sigma(\gamma)-\gamma}+
E^{-\sigma(\gamma)+\gamma},\cr
E^{\gamma}_{-k}a(-1)e^{-\gamma}\in M_2(1)+\langle \sigma(\gamma),a\rangle
e^{\sigma(\gamma)-\gamma}+
\langle \sigma^2(\gamma),a\rangle e^{-\sigma(\gamma)+\gamma}, \quad \mbox{ and} \cr
E^{\gamma}_{-k}\sum_{i=0}^2 \sigma^i(a(-1)e^{-\gamma})
\in M_2(1)+\langle \sigma(\gamma),a\rangle
E^{\sigma(\gamma)-\gamma}+
\langle \sigma^2(\gamma),a\rangle E^{-\sigma(\gamma)+\gamma}. 
\end{array}$$
Therefore, we obtain  
$E^{\sigma(\gamma)-\gamma}, E^{-\sigma(\gamma)+\gamma}\in C_2(V_L^{\sigma})+M_2(1)$.
\prend

For $E^{\mu}$ with 
$\mu=mx+ny$ and $m+n\equiv \pm 1 \pmod{3}$, we need the following lemma. \vspace{-4mm}\\

\begin{lmm}\label{lmm:modulo1}
(1) For $m,n$ with $m+n\equiv 1 \pmod{3}$, there are $\gamma\in L$ satisfying 
$\gamma-\sigma^i(\gamma-\mu)=mx+ny$ for some $i=1,2$ and $\mu\in \{x,y,-x-y\}$ 
such that $\langle \gamma,-\sigma^1(\gamma-\mu)\rangle$ and 
$\langle \gamma,-\sigma^2(\gamma-\mu)\rangle$ are both positive. \\
(2) For $m,n$ with $m+n\equiv 2 \pmod{3}$, there are $\gamma\in L$, $i=1,2$, 
$\mu\in \{-x,-y,+x+y\}$ such that $\gamma-\sigma^i(\gamma-\mu)=mx+ny$, 
$\langle \gamma,-\sigma^1(\gamma-\mu)\rangle>0$ and 
$\langle \gamma,-\sigma^2(\gamma-\mu)\rangle>0$.
\end{lmm}

\pr
We first note that 
for $\gamma=px+qy$ and $-\gamma-x-y$, we have
$$\begin{array}{l}
\langle \sigma(\gamma),-\gamma-x-y\rangle=p^2+q^2-pq+2p-q=(q-\frac{p+1}{2})^2+\frac{3}{4}(p+1)^2-1  \cr
\langle \sigma^2(\gamma),-\gamma-x-y\rangle=p^2+q^2-pq+2q-p=(p-\frac{q+1}{2})^2+\frac{3}{4}(q+1)^2-1,
\end{array}$$
and so the both are positive except $-2\leq p,q\leq 1$. 
For $\mu=mx+ny$ with $m+n\equiv 1\pmod{3}$,  
we may assume $m,n\leq 0$ by taking a conjugate by $<\sigma>$. 
If $\mu=mx+ny\not\in\{x,y,-x-y,-2y\}$, then 
by setting $\gamma=px+qy$ with $q=\frac{-m-n+1}{3}$ and $p=\frac{n-2m+2}{3}$, 
we obtain $\sigma(\gamma)-\gamma-x-y=\mu$ and 
$\langle \sigma(\gamma),-\gamma-x-y\rangle$ and 
$\langle \sigma^2(\gamma),-\gamma-x-y\rangle$ are positive. 
In the case $\mu=-2y$, 
we choose $q=\frac{-m-n+1}{3}$ and $p=\frac{-2n+m+2}{3}$, then 
we have $\sigma^2(\gamma)-\gamma-x-y=\mu$ and 
$\langle \sigma^1(\gamma),-\gamma-x-y\rangle$ and 
$\langle \sigma^2(\gamma),-\gamma-x-y\rangle$ are positive. \\
(2) comes from (1) by replacing $x$ and $y$ by $-x$ and $-y$, 
respectively. 
\prend

By the above lemmas, for any $\mu$, there are $\gamma, \gamma'$ and $k$ such that 
$$
E^{\gamma}_{-2-k}e^{-\gamma'}\in e^{\mu}+e^{\mu'}+M_2(1)e^{\pm x} \quad \mbox{ and so } \quad 
E^{\gamma}_{-2-k}E^{-\gamma'}\in E^{\mu}+E^{\mu'}+M_2(1)E^{\pm x}.$$
We also have 
$$E^{\gamma}_{-2-k+1}\sum_{i=0}^2\sigma^i(a(-1)e^{-\gamma'})
\in \langle a,\gamma\rangle E^{\mu}+\langle a,\sigma(\gamma)\rangle E^{\mu'}+M_2(1)E^{\pm x}, $$
which implies $E^{\mu}, E^{\mu'}\in M_2(1)E^{\pm x}+C_2(V_L^{\sigma})$ for any $\mu$.  
The remaining is to show $a(-n)a'(0)E^{\mu}\in M_2(1)E^{\pm x}+C_2(V_L^{\sigma})$ for 
$n\leq 30$. Actually, we obtain 
$$\begin{array}{l}
E^{\gamma}_{-2-k+1+n}a(-n)a(-n)e^{-\gamma'}\cr
\mbox{}\qquad \in 2n\langle a,\gamma\rangle a(-n)e^{\mu}
+2n\langle a,\sigma(\gamma)\rangle a(-n)e^{\mu'}
+E^{\gamma}_{-2-k+1}e^{-\gamma'}+M_2(1)e^{\pm x}  \qquad \mbox{ and}\cr
E^{\gamma}_{-2-k+1+2n}a(-n)a(-n)a(-n)e^{-\gamma'}\cr
\mbox{}\qquad \in 6n^2\langle a,\gamma\rangle a(-n)e^{\mu}
+6n^2\langle a,\sigma(\gamma)\rangle a(-n)e^{\mu'}
+E^{\gamma}_{-2-k+1}e^{-\gamma'}+M_2(1)e^{\pm x}. 
\end{array}$$
Therefore, we have 
$$\begin{array}{l}
a(-n)e^{\mu}, a(-n)e^{\mu'}\in C_2(V_L^{\sigma})V_L+M_2(1)e^{\pm x} \qquad \mbox{ and so}\cr
a(-n)a'(0)E^{\mu}, a(-n)a'(0)E^{\mu'}\in C_2(V_L^{\sigma})+M_2(1)E^{\pm x} 
\end{array}$$
for $n\leq 30$. This completes the proof of Theorem B.

\section{$\Z_3$-orbifold construction}
Using the known results in \S 2 and Theorem A and B, we will show $\Z_3$-orbifold constructions. 
Let $\Lambda$ be a positive definite even unimodular lattice of rank $N$ 
with a triality automorphism $\sigma$.  We note $8|N$. In this section,  $\xi$ 
denotes $e^{2\pi \sqrt{-1}/3}$. 

Since $\Lambda$ is unimodular, a lattice VOA $V_{\Lambda}$ has exactly one 
simple module $V_{\Lambda}$ and all modules are completely reducible (\cite{D}).  
By \cite{DLiMa},  it has one $\sigma$-twisted module $V_{\Lambda}(\sigma)$ and one 
$\sigma^2$-twisted module $V_{\Lambda}(\sigma^2)$. 
Decompose them as direct sums of simple $V_{\Lambda}^{\sigma}$-modules: 
$$V_{\Lambda}=W^0\oplus W^1\oplus W^2, \quad 
V_{\Lambda}(\sigma)=W^3\oplus W^4\oplus W^5, \quad
V_{\Lambda}(\sigma^2)=W^6\oplus W^7\oplus W^8. $$  

We first show that the weights of elements in $V_{\Lambda}, V(\sigma)$ and 
$V(\sigma^2)$ are in $\Z/3$. 
Set $H=\Lambda^{\sigma}$ and 
$H'=\{u\in \Q H\mid \langle u,h\rangle\in \Z \mbox{ for }h\in H\}$ the dual of $H$. 
Set $s={\rm rank}(H)$, then from the assumption $t=(N-s)/2$ is divisible by three. As it is well known, 
the character $T_{V_H}(\1;\tau)$ of $V_H$ is $\frac{\theta_{H}(\tau)}{\eta(\tau)^s}$, 
where $\eta(\tau)=q^{1/24}\prod_{n=1}^{\infty}(1-q^n)$ is the Dedekind eta-function. 
Since $\Lambda$ is unimodular, $3H'\subseteq H$ and the restriction of $\Lambda$ 
into $\Q H$ covers $H'/H$ and so the weights of elements in $V_H$-modules 
are in $\Z/3$. Hence the powers of $q$ in the character of simple 
$V_H$-modules are all in $-s/24+\Z/3$ and so are those of $q$ 
in $T_{V_H}(\1;-1/\tau)$ by Zhu's theory (\S 2.3.1). 
Since 
$$
T_{V_{\Lambda}}(\sigma,\1;\tau)=q^{-N/24}\frac{\theta_H(\tau)}{\prod_n(1-q^n)^s}\times \frac{1}{\prod_n(1-\xi q^n)^t(1-\xi^{-1}q^n)^t}
=T_{V_H}(1,\1;\tau)\frac{\eta(\tau)^t}{\eta(3\tau)^t}$$
and $T_{V(\sigma)}(1,\1;\tau)$ is a scalar times of 
$$\begin{array}{rl}
T_{V_{\Lambda}}(\sigma,\1;-1/\tau)=&T_{V_H}(1,\1;-1/\tau)\frac{\eta(-1/\tau)^t}{\eta(-3/\tau)^t}
=T_{V_H}(1,\1;-1/\tau)\frac{(\frac{\tau}{\sqrt{-1}})^{t/2}\eta(\tau)^t}{(\frac{\tau}{3\sqrt{-1}})^{t/2}\eta(\tau/3)^t}\cr
=&3^{t/2}T_{V_H}(1,\1;-1/\tau)\frac{\eta(\tau)^t}{\eta(\tau/3)^t}\cr
=&3^{t/2}q^{-2t/24} T_{V_H}(1,\1;-1/\tau)q^{t/9}\frac{\prod_n(1-q^n)^t}{\prod_n(1-q^{n/3})^t}, 
\end{array}\eqno{(5)}$$ 
we have that the powers of $q$ in $T_{V(\sigma)}(1,\1;\tau)$ are in $-N/24+\Z/3$. 
Therefore, we may assume that the weights of elements in $W^i$ are in $i/3+\Z$ for $i\geq 3$. 
In particular, all elements in 
$$\widetilde{V}=W^{0}\oplus W^{3}\oplus W^{6}$$ 
have integer weights. Our next aim is to show that $\widetilde{V}$ has a structure of a 
vertex operator algebra. 
By Theorem B, $V_{\Lambda}^{\sigma}$ is $C_2$-cofinite and all modules are completely reducible 
and $V_{\Lambda}^{\sigma}$ has exactly nine simple modules
$\{W^{i}\mid 0\leq i\leq 8\}$. \vspace{-4mm}\\

\begin{lmm}\label{lmm:allsimplecurrent}
$W^{i}$ are all simple currents, that is, 
$W^{i}\boxtimes W^j$ are simple modules for any $i,j$. 
Moreover, $\widetilde{V}$ is closed by the fusion products. 
\end{lmm}

\pr 
Let determine the entries of the $S$-matrix $(s_{ij})$ of $V_L^{\sigma}$.  
Decompose $S$ into $S=(A_{ij})_{i,j=1,2,3}$ with $3\times 3$-matrices $A_{ij}$. 
Since $S$ is symmetric, $A_{ij}={}^tA_{ji}$.
Simplify the notation, we denote 
$T_{W^i}(v;\tau)$ by $W^i(\tau)$. 
As we explained in \S 2.3.2, there are $\lambda_i\in \C$ $(i=0,1,2)$ 
such that the $S$-transformation shifts 
$$
W^0(\tau) +\xi^i W^1(\tau) +\xi^{2i} W^2(\tau)  
\to \lambda_i(W^{3i}(\tau) +W^{3i+1}(\tau) +W^{3i+2}(\tau) ).$$
Namely, the first three columns of $S$ are 
$$(A_{11}A_{12}A_{13})=\frac{1}{3}\left( \begin{array}{ccccccccc}
\lambda_0&\lambda_0&\lambda_0&\lambda_1&\lambda_1&\lambda_1&\lambda_2&\lambda_2&\lambda_2 \cr
\lambda_0&\lambda_0&\lambda_0&\xi^2 \lambda_1&\xi^2 \lambda_1&\xi^2\lambda_1&\xi \lambda_2&\xi \lambda_2&\xi \lambda_2 \cr
\lambda_0&\lambda_0&\lambda_0&\xi \lambda_1&\xi \lambda_1&\xi \lambda_1&\xi^2\lambda_2&\xi^2\lambda_2&\xi^2\lambda_2 
\end{array}\right).$$ 
Since $S^2$ is a permutation matrix which shifts 
$W$ to its restricted dual $W'$, we get $\lambda_i^2=1$.  
We next consider the characters $\ch(W)=T_{W}(1,\1;\tau)$.  
In this case, since $\ch(W')=\ch(W)$,  we have 
$\ch(W^1)=\ch(W^2)$, $\ch(W^{3+i})$=$\ch(W^{6+i})$ for $i=0,1,2$. 
Clearly, $\{\ch(W^0),\ch(W^1),\ch(W^3),\ch(W^4),\ch(W^5)\}$ is a linearly independent set. 
Since (5) has $q^{1/3+\Z}$-parts, $A_{12}+A_{13}\not=0$ and so $\lambda_1=\lambda_2$.  Similarly, since 
$\ch(W^{3+i})=\ch(W^{6+i})$, we have $A_{22}+A_{23}=A_{32}+A_{33}$.  
Furthermore, since $A_{33}=A_{22}+A_{23}-{}^tA_{23}$ is symmetric, 
$A_{23}$ is symmetric and $A_{22}=A_{33}$. 
As we explained in \S 2.3.2, there are $\mu_i\in \C$ $(i=1,2)$ such that the 
$S$-transformation shifts  
$$W^3(\tau)+\xi^i W^4(\tau)+\xi^{2i} W^5(\tau) 
\to  \mu_i(W^{3i}(\tau)+\xi^2W^{3i+1}(\tau)+\xi W^{3i+2}(\tau)) \qquad \mbox{ for }i=1,2.$$ 
From these information and \S 2.3.2, we know the entries of $S$:
$$ (S_{ij})=
\frac{1}{3}\begin{pmatrix}
\lambda_0&\lambda_0&\lambda_0& \lambda_1&\lambda_1&\lambda_1&       \lambda_1&\lambda_1&\lambda_1 \cr
\lambda_0&\lambda_0&\lambda_0& \xi^2\lambda_1&\xi^2\lambda_1&\xi^2\lambda_1& \xi \lambda_1&\xi \lambda_1&\xi \lambda_1\cr 
\lambda_0&\lambda_0&\lambda_0& \xi \lambda_1&\xi \lambda_1&\xi \lambda_1& \xi^2\lambda_1&\xi^2\lambda_1&\xi^2\lambda_1 \cr 
\lambda_1&\xi^2\lambda_1&\xi \lambda_1& \mu_1&\xi \mu_1& \xi^2 \mu_1 &\mu_2&\xi^2 \mu_2&\xi \mu_2\cr
\lambda_1&\xi^2\lambda_1&\xi \lambda_1& \xi \mu_1&\xi^2 \mu_1&\mu_1&\xi^2 \mu_2&\xi \mu_2&\mu_2  \cr
\lambda_1&\xi^2\lambda_1&\xi \lambda_1& \xi^2 \mu_1&\mu_1&\xi \mu_1&\xi \mu_2&\mu_2&\xi^2 \mu_2  \cr
\lambda_1&\xi \lambda_1&\xi^2\lambda_1& \mu_2&\xi^2 \mu_2&\xi \mu_2&\mu_1&\xi \mu_1&\xi^2 \mu_1  \cr
\lambda_1&\xi \lambda_1&\xi^2\lambda_1& \xi^2 \mu_2&\xi \mu_2&\mu_2&\xi \mu_1&\xi^2 \mu_1&\mu_1  \cr
\lambda_1&\xi \lambda_1&\xi^2\lambda_1& \xi \mu_2&\mu_2&\xi^2 \mu_2&\xi^2 \mu_1&\mu_1&\xi \mu_1,                         
\end{pmatrix}$$
with $\lambda_i^2=\mu_1\mu_2=1$. 
This implies $\overline{S_{ih}}S_{i'h}=1$ and so 
$$ N_{i,i'}^k=\sum_{h} \frac{S_{ih}S_{i'h}S_{hk'}}{S_{0h}}=\sum_h \frac{S_{hk'}}{S_{0h}}.$$
Therefore, $N_{i,i'}^k\not=0$ if and only if $k=0$ and $N_{i,i'}^0=1$. Namely, 
$W^i\boxtimes (W^i)'=V_{\Lambda}^{\sigma}$ for every $i$, If $R\boxtimes W^i$ 
is not simple, then $(R\boxtimes W^i)\boxtimes (W^i)'\cong 
R\boxtimes (W^i\boxtimes (W^i)')\cong R$ is not simple. Therefore, $W^i$ are all simple current.  

By considering the characters, we have:
$$\begin{array}{rl}
T_{V_{\Lambda}}(\sigma,\1;\tau)=&T_{V_H}(1,\1;\tau)\frac{\eta(\tau)^s}{\eta(3\tau)^s}
=\ch(W^0)+\xi\ch(W^1)+\xi^2\ch(W^2)\cr
T_{V_{\Lambda}}(\sigma,\1,-1/\tau)=&\lambda_1\{\ch(W^3)+\ch(W^4)+\ch(W^5)\} \cr
T_{V_{\Lambda}}(\sigma,\1,-1/(\tau+1))=&e^{2\pi \sqrt{-1}N/24}\lambda_1\{\ch(W^3)+\xi\ch(W^4)+\xi^2\ch(W^5)\}\cr
T_{V_{\Lambda}}(\sigma,\1,-1/((-1/\tau)+1))=&e^{2\pi \sqrt{-1}N/24}\lambda_1\mu_1\{\ch(W^3)+\xi^2\ch(W^4)+\xi\ch(W^5)\}
\end{array}$$
from the above $S$-matrix. On the other hand, since 
$$\begin{array}{rl}
T_{V_{\Lambda}}(\sigma,\1,-1/((-1/\tau)+1))=&T_{V_{\Lambda}}(\sigma,\1,-1-\frac{1}{\tau-1}) \cr
=&e^{-2\pi \sqrt{-1}N/24}T_{V_{\Lambda}}(\sigma,\1,-1/(\tau-1)) \cr
=&e^{-4\pi \sqrt{-1}N/24}\lambda_1\{\ch(W^3)+\xi^2\ch(W^4)+\xi \ch(W^5)\}
\end{array}$$
we have $\mu_1=e^{-6\pi \sqrt{-1}N/24}=1$ and $\mu_1=1$ since $8|N$.  
Then the $S$-matrix implies $\lambda_1=\lambda_0$ and 
$W^3\boxtimes W^3=W^6$ and $W^3\boxtimes W^6=W^0$.

Therefore, 
$$ \widetilde{V}=W^{0}\oplus W^{3}\oplus W^{6}$$ 
is a direct sum of simple current $V$-modules $W^{3i}$ with integer weights and 
$W^{3i}\boxtimes W^{3j}=W^{3k}$ with $i+j\equiv k \pmod{3}$. 
In order to prove that $\widetilde{V}$ has a VOA-structure,  
we will prove a more general statement.  

\begin{prn} Let $V$ be a $C_2$-cofinite VOA of CFT-type and all $V$-modules are 
completely reducible. Let $W=\oplus_{i=0}^n W^i$ be a direct sum 
of simple $V$-module $W^i$ with integer weights and we assume 
$W^i\boxtimes W^j=W^{k}$ for $i+j \equiv k \pmod{n}$, $W^0=V$ and $n$ is odd. 
Then $W$ has a VOA structure. 
\end{prn}

\pr
Simplify the notation, $W^s$ denotes $W^i$ for $s\in \Z$ with $s\equiv i \pmod{n}$. 
Let $\CY^{i,j}$ be a nonzero intertwining operator of type $\binom{W^{i+j}}{W^i\quad W^j}$.  
We choose $\CY^{2,i}$ so that 
$$E(\langle d', \CY^{1,i+1}(w,x)\CY^{1,i}(u,y)a\rangle) 
=E(\langle d', \CY^{2,i}(\CY^{1,1}(w,x-y)u,y)a\rangle)$$
for any $a\in W^i$, $w,u\in W^1$ and $d'\in (W^{i+2})'$.
Set an intertwining operator $\CY$ of type $\binom{W}{W^1 \quad W}$ by 
$$
\CY(w,z)=
\begin{pmatrix}0& \cdots&0&\CY^{1,n-1}(w,z) \cr
\CY^{1,0}(w,z)&0\cdots&0&0 \cr
\vdots &\ddots&\vdots&\vdots \cr
0&\cdots 0&\CY^{1,n-2}(w,z)&0\end{pmatrix},$$
where $w\in W^1$. We note that a vertex operator $Y^W$ of $V$ on $W$ is given by 
$$Y(v,z)=\begin{pmatrix} Y^{0,0}(v,z)&\cdots&0 \cr 
\vdots&\ddots & \vdots\cr
0&\cdots & Y^{0,n-1}(v,z) \end{pmatrix}.$$ 
By the definition, it is easy to see that 
$\CY(w,z)$ satisfies Commutativity with $Y^W(v,z)$ for any $w\in W$ and $v\in V$. 
It is also easy to check $\CY(L(-1)w,z)=\frac{d}{dz}\CY(w,z)$. 

Our next aim is to prove that $\CY(w,z)$ satisfies Commutativity with itself.
We note that $\CY^{i,j}$ are all integer power series. Therefore, 
it is sufficient to show 
$$ E(\langle d', \CY^{1,i+1}(w,x)\CY^{1,i}(u,y)a \rangle)
=E(\langle d', \CY^{1,i+1}(u,y)\CY^{1,i}(w,x)a \rangle)$$
for any $i=0,\cdots,n-1$, $d'\in (W^{i+2})'$, $a\in W^i$, and $w,u\in W^1$.
 
Recall a skew-symmetric intertwining operator 
$\sigma_{12}(\CY^{1,1})(w,z)u=e^{L(-1)z}\CY^{1,1}(u,-z)w$. 
Since $\dim \CI_{W^1,W^1}^{W^2}=1$ and all $W^i$ have integer weights, 
we have $\sigma_{12}^2=1$ on $\CI_{W^1,W^1}^{W^2}$ and so there is 
$\lambda\in \{\pm 1\}$ such that $\sigma_{12}(\CY^{1,1})=\lambda \CY^{1,1}$. 
Therefore we have:
$$\begin{array}{rl}
E(\langle d', \CY^{1,i+1}(w,x)\CY^{1,i}(u,y)a\rangle) 
=&E(\langle d', \CY^{2,i}(\CY^{1,1}(w,x-y)u,y)a\rangle) \cr
=&E(\langle d', \CY^{2,i}(e^{L(-1)(x-y)}\sigma_{12}(\CY^{1,1})(u,y-x)w,y)a\rangle) \cr
=&E(\langle d', \CY^{2,i}(\sigma_{12}(\CY^{1,1})(u,y-x)w,x)a\rangle) \cr
=&E(\langle d', \lambda\CY^{2,i}(\CY^{1,1}(u,y-x)w,x)a\rangle) \cr
=&E(\langle d', \lambda\CY^{1,i+1}(u,y)\CY^{1,i}(w,x)a\rangle).
\end{array}$$
Irritating it $n$ times, we obtain   
$$\begin{array}{l}
E(\langle e', \CY^{1,n-1}(d^1,x_1)\cdots \CY^{1,0}(d^n,x_n)\CY^{1,n-1}(a,z_1)\cdots \CY^{1,0}(c,z_{n-1})v\rangle) \cr
=\lambda^{n^2}E(\langle e', \CY^{1,0}(a,z_1)\cdots \CY^{1,1}(c,z_{n-1})
\CY^{1,n-1}(d^1,x_1)\cdots \CY^{1,0}(d^n,x_n)v\rangle)
\end{array}$$
for $e'\in (W^1)'$, $a,\ldots,c,d^i\in W^1$ and $v\in V$. On the other hand, from the associativity 
$\overbrace{W^1\boxtimes \cdots \boxtimes W^1}^n=V$, there is $0\not=\mu\in \C$ 
such that 
$$\begin{array}{l}
E(\langle e', \CY^{1,n-1}(d^1,x_1)\cdots \CY^{1,0}(d^n,x_n)\CY^{1,n-1}(a,z_1)\cdots \CY^{1,0}(c,z_{n-1})v\rangle)\cr
\mbox{}\qquad =E(\langle e', \CY^{1,n-1}(d^1,x_1)\cdots \CY^{1,0}(d^n,x_n)\times \cr
\mbox{}\qquad \qquad \times \mu Y^V(\CY^{n-1,1}\cdots\CY^{1,1}(a,z_1-z_2)b,
\cdots (z_{n-2}-z_{n-1})c,z_{n-1})v\rangle).\end{array}$$
for any $v\in V$. Since $\CY^{i,j}$ satisfies Commutativity with $Y^W$, we have 
$$\begin{array}{l}
E(\langle e', \CY^{1,n-1}(d^1,x_1)\cdots \CY^{1,0}(d^n,x_n)\CY^{1,n-1}(a,z_1)\cdots \CY^{1,0}(c,z_{n-1})v\rangle) \cr
\mbox{}\qquad =E(\langle e', \CY^{1,n-1}(d^1,x_1)\cdots \CY^{1,0}(d^n,x_n)\times \cr
\mbox{}\qquad\qquad \times \mu Y^V(\CY^{n-1,1}\cdots \CY^{1,1}(a,z_1-z_2)b,
\cdots (z_{n-2}-z_{n-1})c,z_{n-1})v\rangle)\cr
\mbox{}\qquad =E(\langle e', \mu Y^V(\CY^{n-1,1}\cdots \CY^{1,1}(a,z_1-z_2)b,
\cdots (z_{n-2}-z_{n-1})c,z_{n-1})\times \cr
\mbox{}\qquad\qquad \times \CY^{1,n-1}(d^1,x_1)\cdots \CY^{1,0}(d^n,x_n)v\rangle)\cr
=E(\langle e', \CY^{1,0}(a,z_1)\cdots \CY^{1,1}(c,z_{n-1})
\CY^{1,n-1}(d^1,x_1)\cdots \CY^{1,0}(d^n,x_n)v\rangle), 
\end{array}$$
which implies $\lambda^{n^2}=1$.  Since $n$ is odd, we have $\lambda=1$.
Thus, $\CY$ satisfies Commutativity with itself. 
By using the normal products, $\CY$ and $Y^W$ generate a local system $\CO$ 
with a Virasoro element $Y^W(\omega,z)$.  
Since $V$ is a subVOA of $\CO$ and its modules are completely reducible, 
$\CO$ is a direct sum $\CO=\oplus T^j$ of simple $V$-modules. 
Clearly, the action of $\CO$ on $W$ induces intertwining operators of type 
$\binom{W}{\CO\quad W}$. Since $V\boxtimes V=V$, 
the multiplicity of a $V$-module $V$ in $\CO$ is one and so we have 
$\CO\cong W$ as $V$-modules. Therefore, $W$ has a VOA structure. \prend

\subsection{The moonshine VOA}
Let's apply the above construction to the Leech lattice $\Lambda$ and a fixed point 
free automorphism $\sigma$ of $\Lambda$ of order three. 
Then a trace function $T_{V_{\Lambda}}(\sigma,\1;\tau)$ of $\sigma$ on $V_{\Lambda}$ is 
$$
q^{-1}(\frac{1}{\prod_{n=1}^{\infty}(1-\xi q^n)})^{12}
(\frac{1}{\prod_{n=1}^{\infty}(1-\xi^2 q^n)})^{12}=q^{-1}(\frac{1}{\prod_{n=1}^{\infty}(1+q^n+q^{2n})})^{12}=\frac{\eta(\tau)^{12}}{\eta(3\tau)^{12}}.
$$
Hence, a character function of the twisted module $V_{\Lambda}(\sigma)$ is 
$$ \ch(V_{\Lambda}(\sigma))=\ch(W^{3})+\ch(W^{4})+\ch(W^{5})=T_{V_{\Lambda}}(\sigma,\1,-1/\tau)=
3^6q^{-1}q^{4/3}\frac{\prod_{n=1}^{\infty}(1-q^n)^{12}}{\prod_{n=1}^{\infty}(1-q^{n/3})^{12}},$$
which implies that $W^{3}$ (also $W^{6}$) has no elements of weight $1$ and 
$\ch(\widetilde{V}_{\Lambda})=J(\tau)$. 

By an easy calculation, $\dim W^3_2=3^6(12+12+\binom{12}{2})=65610$ and so a triality automorphism of 
$\widetilde{V}$ defined by $e^{2\pi \sqrt{-1}i/3}$ on $W^{3i}$ for $i=0,1,2$ is corresponding to 
${\rm 3B}\in {\mathbb M}$ if 
$\widetilde{V}\cong V^{\natural}$. 

\subsection{A new VOA No.32 in Schellekens' list}
We next start from a Niemeier lattice $N$ of type $E_6^4$ and 
a triality automorphism $\sigma$ which 
acts on the first entry $E_6$ as fixed point free and 
permutes the last three $E_6$, where we choose $<(0,1,1,1),(1,1,2,0)>$ 
as a set of glue vectors of $N$ for $E_6^4$. We note that 
since $E_6$ contains a full sublattice $A_2^3$, 
$E_6$ has a fixed point free automorphism of order three. 
Since $t=9$, in order to determine the dimension of $(\tilde{V}_N)_1$, 
it is enough to see the constant term of 
$q^{6/24}\frac{\Theta_H(-1/\tau)}{\eta(-1/\tau)^6}$. 
Since the fixed point sublattice $H$ is isomorphic to $\sqrt{3}E_6^{\ast}$, 
we have  
$$\Theta_H(\tau)
=\frac{1}{3}[\phi_0(\tau)^3+\frac{1}{4}\{3\phi_0(3\tau)-\phi_0(\tau)\}^3], $$
where $\phi_0(\tau)
=\theta_2(2\tau)\theta_2(6\tau)+\theta_3(2\tau)\theta_3(6\tau)$ 
and $\theta_2(\tau)=\sum_{m\in \Z}q^{(m+1/2)^2}$ and 
$\theta_3(\tau)=\sum_{m\in \Z}q^{m^2}$, see \cite{CS}. 
Applying $\phi_0(-1/\tau)
=\frac{\tau}{i\sqrt{3}}\phi_0(\tau/3)$,  we have 
$$\frac{\Theta_H(-1/\tau)}{\eta(-1/\tau)^6}=\frac{1}{9\sqrt{3}}q^{-1/4}+\cdots $$
and so 
$$\dim (\tilde{V}_N)_1=(6\times 12)/3+6+6\times 12+2\times\{3^{9/2}3^{-5/2}\}=120.$$ 
Clearly, from the construction we know that 
$(\tilde{V}_N)_1$ contains $A_2^3E_{6,3}$ as a subring.  
Therefore, $\tilde{V}_N$ is a new 
vertex operator algebra No 32 in the list of Schellekens \cite{S}. \\

\noindent
{\bf Acknowledgement}  
The author would like to thank K.~Tanabe for his right questions.

\end{document}